\documentclass[12pt,a4paper,twoside,final,notitlepage, leqno]{article}
\usepackage[english]{babel}
\usepackage[T1]{fontenc}
\usepackage{epsfig, graphicx, amssymb}
\usepackage{amsmath,amsthm,epsfig,amsfonts}
\usepackage{float}
\usepackage{color}
\def\br {\break}

\newcommand{\moneq}{\vspace*{-7pt} \begin{equation} \displaystyle }
\newcommand{\moneqstar}{\vspace*{-6pt} \begin{equation*} \displaystyle }
\newcommand{\monendstar}{\vspace*{-6pt} \end{equation*}   }
\newcommand{\monend}{\vspace*{-7pt} \end{equation}   }
\newcommand{\moneqarraystar}{ \begin{eqnarray*} \displaystyle }
\newcommand{\monendarraystar}{ \end{eqnarray*}   }

\newcommand{\dd}{{\rm d}}

\newcommand{\RR}[0]{\mathbb{R}}

\definecolor{vertfonce}{rgb}{0.0, 0.5, 0.0}

\def\section*#1{}

\usepackage{fancyhdr}
\renewcommand{\headrulewidth}{0pt}

\begin{document}

\fancypagestyle{plain}{ \fancyfoot{} \renewcommand{\footrulewidth}{0pt}}
\fancypagestyle{plain}{ \fancyhead{} \renewcommand{\headrulewidth}{0pt}}

~

  \vskip 2.1 cm

\centerline {\bf \LARGE Numerical approximations}

\bigskip

\centerline {\bf \LARGE of a lattice Boltzmann scheme}

\bigskip

\centerline {\bf \LARGE with a family of partial differential equations}

 \bigskip  \bigskip \bigskip

\centerline { \large   Bruce M. Boghosian$^{ab}$,   Fran\c{c}ois Dubois$^{cd}$ and Pierre Lallemand$^{e}$}

\smallskip  \bigskip

\centerline { \it  \small
  $^a$ Department of Mathematics, Tufts University, Medford, MA, 02155, USA.}

\centerline { \it  \small
  $^b$ present address: American University of Armenia, }

\centerline { \it  \small 40 Baghramyan Avenue, Yerevan 0019, Armenia.}

\centerline { \it  \small
  $^c$   Laboratoire de Math\'ematiques d'Orsay, Facult\'e des Sciences d'Orsay,}

\centerline { \it  \small   Universit\'e Paris-Saclay, France.}

\centerline { \it  \small
$^d$    Conservatoire National des Arts et M\'etiers, LMSSC laboratory,  Paris, France.}


\centerline { \it  \small
 $^e$ Beijing Computational Science Research Center, Haidian District, Beijing 100094,  China.}


\bigskip  \bigskip

\centerline {02 August 2024 
  {\footnote {\rm  \small $\,$ Computers and Fluids (Mesoscopic methods and their applications to CFD), volume 284,
    article~106410~(14 pages), november 2024. 
This contribution has been presented 
at the 19th International Conference for Mesoscopic Methods in Engineering and Science,
Mount Qing-Cheng (Chengdu, Sichuan, China) the 25 July 2023
and at Institut Henri Poincar\'e the 04 October 2023.}}}




\bigskip  \bigskip
\noindent {\bf \large Abstract}

\noindent
In  this contribution, we address 
the  numerical solutions of high-order
asymptotic equivalent     partial differential equations 
with the results of a  lattice Boltzmann scheme for an  
inhomogeneous advection problem in one spatial dimension.
    We  first derive a family of
equivalent partial differential equations
at various orders, 
and we compare the lattice Boltzmann experimental
results with  
a spectral approximation of the differential equations. 
For an unsteady situation, we show that the initialization scheme at a sufficiently  
high order
of the microscopic moments plays
a crucial role to observe an asymptotic error consistent   
with the order of approximation.
For a stationary long-time  
limit, we observe that the measured asymptotic error converges   
  with a reduced order of precision compared to the one suggested by asymptotic analysis.

\noindent

\newpage

\noindent {\bf \large    1) \quad  Introduction} 

\fancyhead[EC]{\sc{Bruce M. Boghosian,  Fran\c{c}ois Dubois and Pierre Lallemand}}
\fancyhead[OC]{\sc{Numerical approximations of a lattice Boltzmann scheme}}
\fancyfoot[C]{\oldstylenums{\thepage}}

\smallskip \noindent
The classical framework for numerical simulation 
starts from partial differential equations.
After discretization with some numerical method (finite differences, finite elements, {\it etc.}),
numerical software is developed. Then  an approximate solution
of the original partial differential equation is computed.
A large number of high-quality books  exist on this subject.
We refer to the works of 
Oden  and Reddy \cite{OR76}, 
Ferziger and Peri\'c \cite{FP96},   
Lucquin and  Pironneau~\cite{LP98}, among others. 

\smallskip \noindent
With cellular automata and lattice Boltzmann schemes, this paradigm is reversed.
The computing algorithm is the starting point of the study. Then an asymptotic analysis
is conducted to derive the underlying continuous equations.
The reader can consult, {\it e.g.}, the books of
Rothman and  Zaleski \cite{RZ97},
Succi \cite{Su01},
Guo and Shu \cite{GS13},
or  Kr\"uger {\it et al.}~\cite{KKKSSV17}. 

\smallskip \noindent
Our approach is to derive  a physical model  from the algorithm.
In this work, we adopt the paradigm  of  multirelaxation lattice Boltzmann
schemes~\cite {DDH92}.
Other variants  
are possible and  there are rational ways to
proceed from a partial differential equation system with suitable structure to a kinetic formulation,
and hence to a lattice Boltzmann scheme by further
discretisation in space and time, e.g. by integration along characteristics
(He and Luo \cite{HL97}, He {\it et al.} \cite{HSD98},  Dellar \cite{De13}).
The  classical approach for approaching the continuous Boltzmann equation is 
the Chapman-Enskog method \cite{CC39}. It has been  revisited in  Chen and Doolen \cite{CD98}
and Qian and Zhou \cite{QZ00} to take into
consideration the discrete aspects of space and time with cellular automata and lattice Boltzmann schemes.

\smallskip \noindent
  One important remark has to do with  
  the choice of  scaling. In this contribution,
  we suppose  acoustic scaling: the ratio between the spatial step 
  and the time step is fixed.
  Then it is possible to derive asymptotic partial differential equations in terms
  of a purely numerical parameter, {\it e.g.} the spatial step  
  to fix the ideas.
Taylor expansions allow us  
to derive equivalent partial differential equations \cite{Du08,fd09}.
Dubois {\it et al.} 
have  established that this Taylor expansion method is equivalent to the Chapman-Enskog approach~\cite{DBL23}.
 This asymptotic expansion is obtained by formal arguments
 and can be compared to the truncation error of a finite difference scheme.
Observe that the relaxation coefficients are supposed fixed when we adopt the acoustic scaling hypothesis.

\smallskip \noindent
When an asymptotic partial differential equation is known, it is possible to fit some
parameters of the scheme to obtain super convergence.
This was done by d'Humi\`eres and Ginzburg~\cite {HG09},
Augier {\it et al.} \cite{ADGG13}, Dubois and Lallemand   \cite{DL09,DL11}
and  Otomo {\it et al.} \cite{OBD17}.
%
In the present contribution, 
we fix the relaxation coefficients and consider the spatial step 
as tending to zero.

\smallskip \noindent
 An asymptotic expansion is not identical to a convergence result.
 With the acoustic scaling  at second order to fix  ideas,  
 the diffusivity is proportional to the spatial step 
 and  vanishes  as  the spatial step  
 tends to zero.
 Moreover,  the truncation error is {\it a priori} not  identical to the mathematical error.
 Some ad hoc stability and consistency are necessary to establish convergence, 
as classically established by Lax and Richtmyer in \cite{LR56}. 
After the pioneering work of Dellacherie \cite{De14},
Boghosian {\it et al.}  have established in  \cite{BDGLT18a,BDGLT18b}
that the lattice Boltzmann method does not converge 
for a simple heat equation.
The acoustic scaling necessitates modification of   
the relaxation coefficients to maintain fixed  
diffusivity. In consequence, the hypothesis of fixed relaxation coefficients
is no longer valid,  
and the target partial differential equation is not valid for very small spatial steps.   
 Nevertheless, the asymptotic expansion with  acoustic scaling is correct.
 It is simply defining a partial differential equation that  mimics 
 the lattice  Boltzmann scheme at a specific order of accuracy 
 when the relaxation coefficients are fixed and the spatial step  
 tends  
 to zero.

\smallskip \noindent
Our methodology to develop an asymptotic analysis is based on the ABCD method of analysis~\cite{Du22}
  enforced by the equivaence with the Chapman-Enskog expansion \cite{DBL23}.
  It allows a fourth-order 
  analysis in a wide class of nonlinear  lattice Boltzmann schemes.
Thus an important question is the comparison between the simulation with a lattice Boltzmann scheme
and reference solutions of the equivalent partial differential equations.
The objective of this contribution  is to addres 
the  numerical solution  of high-order 
equivalent  asymptotic   partial differential equations 
with a  lattice Boltzmann scheme.
We work with an elementary D1Q3 one-dimensional lattice Boltzmann scheme
with a prescribed sinusoidal velocity.
We treat the reference asymptotic partial differential equations
  with an approach  reminiscent of early kinematic
dynamo simulations from the 1980s using pure spectral methods to solve the magnetic induction equation
in ``ABC'' ﬂow ({\it e.g.} Galloway and Frisch \cite{GH86}).
With this    spectral approach we can numerically solve  
  with great precision the family of equivalent
partial differential equations at various orders.

As a classical tool of analysis, we use the one-point  
spectral analysis developed in
Lallemand and Luo \cite{LL00} (see also  Simonis {\it et al.} \cite{SHKDK21}). 
We use also Arnoldi iterations \cite{Ar51} for the determination of global modes
for an entire mesh, as   Verberg and Ladd \cite{VL99} to  compute
steady solutions of linear lattice Boltzmann schemes for Stokes ﬂow using a Krylov-space method
of Leriche {\it et al.} \cite {LLL08} 
for the determination of Stokes eigenmodes  in a cubic domain. 

For many simulations, 
the initialisation   of the lattice Boltzmann state is taken to be an  
equilibrium.
This simple approach has been enriched by first-order   
initialization  initially  suggested to  our knowledge by Mei {\it et al.}  \cite{MLLH06}. 
Within   
the framework of Bellotti {\it et al.}  \cite{BGM22},
a lattice Boltzmann scheme is revisited as a multistep method
for the conserved variables. The equivalent partial differential equations have
been established at second-order accuracy with this framework under acoustic scaling~\cite{Be23}.
Last but not least, the  contribution of Mei {\it et al.} \cite{MLLH06} 
at first-order accuracy has been revisited by Bellotti ~\cite{Be24}.

\smallskip \noindent
The outline  for this work is as follows.
In Section~2, we study the reference model: the  advection equation  in one spatial dimension
with a given cosine velocity field. The method of characteristics yields an analytic solution.
In Section~3, we present our variant of the  D1Q3 lattice Boltzmann scheme,
introduced initially by Broadwell \cite{Br64} 
in the  context of  simple discrete-velocity gases.
In the lattice Boltzmann framework,  dynamics is captured with particles
and the relaxation process occurs in the space of moments~\cite {DDH92}.
They are  divided into two families: the conserved moments and the microscopic variables
in the denomination proposed by Gatignol~\cite{Ga87}.
%

\smallskip \noindent
Then,  in Section 4 we present the ``ABCD'' asymptotic analysis  \cite{Du22, DBL23}.
From the precise algebraic expression  of a multirelaxation lattice Boltzmann
scheme~\cite {DDH92,HGKLL02}, we  derive from a formal exponential expression
a set of  equivalent partial differential equations up to fourth
order accuracy. Here we adapt the underlying algebra first to the case of an  
inhomogeneous advection problem in one spatial dimension, and we  derive a family of
equivalent partial differential equations
at various orders.

\smallskip \noindent
In Section~5, the Fourier series method is adapted to treat in a precise way the
case of a cosine advective field.
Then the unsteady evolution is presented in Section~6. A first result is relative to  a constant velocity
and an  initial sinusoidal wave. Then we take into account a cosine advection velocity with a
sinusoidal or a constant initial condition.
We compare the lattice Boltzmann results and a spectral approximation of the differential equations. 
An interesting phenomenon of lack of convergence is encountered. 
This motivates the next Section relative to the initialization of microscopic moments.
In Section 8, we present our numerical experiments with a detailed asymptotic analysis.
Various parameters are considered: the type of problem, with constant or cosine advective velocity,
the approximation order of the partial differential equation, the number of mesh points and
the initialization process.
In Section 9, we study the long-time asymptotics.   
We observe that the measured asymptotic error is still converging.

\smallskip \noindent
This work is the result of  conversations  in Medford (MA, USA) during summer 2018, and then in Paris in spring 2019.
Independent numerical experiments were done during COVID in spring  2020,
and complementary work in Beijing in summer  2023.

\bigskip \bigskip    \noindent {\bf \large    2) \quad  Advection with harmonic velocity in one space dimension} 

\smallskip \noindent
We introduce a reference length $ \, L > 0 $,  a final time $\, T \, $
and a reference scale velocity $ \,  \lambda  > 0 $.
For a given scalar $ \, U \in \RR \, $ and for $ \, 0 \leq x \leq L $, we consider the regular periodic  velocity field
\moneq \label{u-cosinus}
u(x) =  \lambda \, U \, \cos(k\, x) \,,\,\,  k = {{2 \pi}\over{L}} .
\monend 
%
%
The linear inhomogeneous advection equation is the  first-order partial differential equation 
\moneq \label{advection-vitesse-variable} 
  {{\partial \rho}\over{\partial t}} + \lambda \, {{\partial}\over{\partial x}} \big[ U \, \cos(k\, x) \, \rho \big] = 0 .
\monend 
We introduce a periodic function $ \, [0,L] \ni x \longmapsto \rho_0(x) \in \RR \,  \, $ as  an initial condition
\moneq \label{advection-condition-initiale}
\rho(x,\, 0)  = \rho_0(x) . 
\monend
Moreover, we suppose periodic boundary conditions throughout this study.

\bigskip \noindent
{\bf Proposition 1. Method of characteristics}

\noindent 
The differential equation associated with the method of characteristics for the partial
differential equation (\ref{advection-vitesse-variable}) is written 
\moneq \label{advection-vitesse-variable-edo} 
    {{\dd X}\over{\dd t}} = \lambda \,u(X(t)) \equiv \lambda \,U \, \cos \big(k\, X(t) \big) .
\monend
With the initial condition  $ \, X(0) = x_0 \, $ with  $ \, 0 \leq x_0 \leq L $,  
the solution is:
\moneq \label{solution-caracteristiques}
{\rm cotg} {{\pi \, X}\over{L}} =  \displaystyle {{{\rm th}{{\pi \, t}\over{T}} + {\rm cotg} {{\pi \, x_0}\over{L}}  }\over
{1 + {\rm th}{{\pi \, t}\over{T}} \, {\rm cotg} {{\pi \, x_0}\over{L}}     }}
\monend 
with $\,\, \lambda \, U \equiv {{L}\over{T}} $, 
$ \,\, \displaystyle {\rm th} \varphi \equiv
{{\exp \varphi - \exp (-\varphi)}\over{\exp \varphi + \exp (-\varphi)}} \, $ and
$ \,  {\rm cotg} \varphi \equiv {{1}\over {\displaystyle{\rm tan} \varphi}} $. 

\bigskip \noindent
{Proof of Proposition 1.}

\noindent 
If $ \, t = 0 $, then $ \,   {\rm cotg} {{\pi \, X}\over{L}} = {\rm cotg} {{\pi \, x_0}\over{L}} \, $
and $ \, \pi \, {{X-x_0}\over{L}} \, $ is a multiple of $ \, \pi $. Then 
the position $ \, X = x_0 \, $ is well defined in the interval $ \, [0 ,\, L ] $.
Moreover, we have the following calculation:
\moneqstar  -{{1}\over{{{\rm sin}^2} {{\pi \, X}\over{L}}}} \, {{\pi}\over{L}} \, {{\dd X}\over{\dd t}} =  {{\pi}\over{T}} \,
{{(1-{\rm th}^2{{\pi \, t}\over{T}})\,  (1- {\rm cotg}^2 {{\pi \, x_0}\over{L}})}
  \over{(1 + {\rm th}{{\pi \, t}\over{T}} \, {\rm cotg} {{\pi \, x_0}\over{L}})^2 }}
= - {{\pi}\over{T}} \, {{1}\over{{{\rm sin}^2} {{\pi \, X}\over{L}}}} \, \cos  {{2 \, \pi \, X}\over{L}}
\monendstar
Then 
$ \,\,  {{\dd X}\over{\dd t}} =  {{L}\over{T}} \,  \cos  {{2 \, \pi \, X}\over{L}} \,\, $
and the differential equation  (\ref{advection-vitesse-variable-edo}) is satisfied.
\hfill $\square $


\bigskip 
\noindent
{\bf Proposition 2. Algebraic solution of the inhomogeneous advection equation}

\noindent
Given $ \, x \in   [0  ,\, L ] \,\, $ and $ \, t > 0 $,
the solution $ \, \rho(x,\, t) \, $ of the equation (\ref{advection-vitesse-variable}) satisfying
the initial condition~(\ref{advection-condition-initiale})
is given by the relation
\moneqstar
\rho(x,\, t) \, \cos \Big({{2 \pi x}\over{L}} \Big) = \rho_0(x_0) \, \cos \Big({{2 \pi x_0}\over{L}} \Big) 
\monendstar
where $ \, x_0 \, $ satisfies 
\moneq \label{pied-caracteristique}
{\rm cotg} {{\pi \, x_0}\over{L}} 
= {{ {\rm cotg} {{\pi \, x}\over{L}} - {\rm th}{{\pi \, t}\over{T}}  }\over
{1 - {\rm th}{{\pi \, t}\over{T}} \,\, {\rm cotg} {{\pi \, x_0}\over{L}}     }} .
\monend 

\bigskip \noindent
{Proof of Proposition 2.}

\noindent 
We search
for  
the foot $ \, x_0 \, $  of the characteristic $ \, t \longmapsto X(t) \, $ (\ref{advection-vitesse-variable-edo}) 
such that $ \, X(0) = x_0 \, $ and $ \, X(t) = x $.
The characteristic passing through point $ \, x \, $  at time $ \, t \, $ 
  then satisfies the two relations
$\,   X(0) = x_0 \, $ and $ \, \, X(t) = x $.
First, we deduce from the partial differential equation  (\ref{advection-vitesse-variable}) that the product
$ \, \rho(x,\, t) \, \cos \big({{2 \pi x}\over{L}} \big) \, $ remains constant.
Second, 
from the relation~(\ref{solution-caracteristiques}),   we deduce the relation~(\ref{pied-caracteristique})
for defining $ \, x_0 $. 
\hfill $\square $

\renewcommand{\thefigure}{1}
\begin{figure}    [H]  \centering
\centerline  {\includegraphics[width=.28\textwidth]   {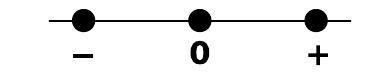}}
\caption{D1Q3 lattice Boltzmann scheme}
\label{fig-01} \end{figure}

\bigskip \bigskip
\noindent {\bf \large    3) \quad  D1Q3 lattice Boltzmann scheme} 

\smallskip \noindent   
The scale velocity $ \, \lambda  > 0  \, $ is now equal to the ratio between the spatial step  
$ \, \Delta x \, $ and the time step $ \, \Delta t $:
\moneqstar
\lambda = {{\Delta x}\over{\Delta t}} .
\monendstar
From the particle distribution $ \, f \equiv \big(f_+ ,\,   f_0 ,\,  f_- \big)^{\rm t} \, $
presented in  Figure \ref{fig-01}, 
we construct a single 
conserved moment $ \, W \, $ denoted ``density'' in the following:
$ \, \rho = f_+ + f_0 + f_- $. 
We have also two non-conserved microscopic moments  $ \, Y = \big( J ,\, e \big)^{\rm t}  \, $
with the  ``momentum''   $ \, J =  \lambda \,  f_+  -\lambda \, f_-  \, $
and the  ``energy''   $ \, e =  \lambda^2 \, ( f_+   - 2 \, f_0 + f_-  ) $.
Then the family of moments
$ \, m \equiv (W \,,\,\, Y) \, $ is linked to the particles $ \, f \, $ with the 
d'Humi\`eres \cite{DDH92} matrix $ \, M  $: $ \, m \equiv M \, f $, 
with
\moneqstar
M = \begin{pmatrix}   1 & 1 & 1 \\  \lambda & 0 & - \lambda \\ \lambda^2 & -2 \, \lambda^2 & \lambda^2  \end{pmatrix} .
\monendstar 
%
%
%
For an inhomogeneous linear equilibrium $ \, Y^{\rm eq} = \Phi(W) = E(x) \, W $,
the equilibrium matrix~$ \, E(x) \, $ is a function of space. In the case of an advective field $ \, u(x) \, $
proposed in the relation (\ref{u-cosinus}), we have 
\moneqstar
E(x) = \begin{pmatrix} \lambda \, U \, \cos(k\,x) \\ \lambda^2 \, \alpha \end{pmatrix}
\monendstar
with a coefficient $ \, \alpha =  -1 \, $ in our numerical experiments. 

\smallskip \noindent
The  relaxation $\, Y \longmapsto Y^* \, $ of the nonconserved moments $ \, Y \, $ is classical:
\moneqstar  \left\{ \begin{array} {l}
J^* = J + s \, (J^{\rm eq} - J) = (1-s)\, J + s \, \lambda \, U \, \cos(k\,x) \, \rho \\ 
e^* = e + s' \, (e^{\rm eq} - e) = (1-s')\, e + s' \, \lambda^2 \, \alpha  \, \rho
\end{array} \right. \monendstar 
and we have chosen
$ \, \, s = 1.5 $, $\, s' = 1.2 \,\, $ in our reference numerical experiments.
Observe that the parameters 
$ \, U $, $ \, \alpha $, $ \, s \, $ and $ \, s' \, $ are without dimension.
We set finally $ \, m^* = ( \rho \,,\,\, J^*  \,,\,\, e^*  )^{\rm t} $. 

\smallskip \noindent
The collision step is defined according to
$ \, f^* =  M^{-1}  \, m^*  $,  and the exact propagation of particles along the
characteristic directions $ \, \lambda \,,\,\, 0 \,,\,\, -\lambda \, $
of the D1Q3 scheme:
\moneqstar  \left\{ \begin{array} {l}
f_+ (x,\, t+\Delta t) = f_+^*(x-\Delta x,\, t)  \\ 
f_0(x,\, t+\Delta t) \, = f_0^*(x ,\, t )  \\ 
f_-(x,\, t+\Delta t) = f_-^*(x+\Delta x,\, t ) 
\end{array} \right. \monendstar
is well known (see {\it e.g.} \cite{DDH92}). The solution of this 
lattice Boltzmann scheme can be approached by  the  first order   equivalent partial differential equation 
\moneqstar
{{\partial \rho}\over{\partial t}} + \lambda \, {{\partial}\over{\partial x}} \big[ U \, \cos(k\, x)  \, \rho \big] = {\rm O}(\Delta x) .
\monendstar 
Therefore, it is natural to compare the numerical solution of the D1Q3 lattice Boltzmann scheme with the
exact solution of the   inhomogeneous  advection equation (\ref{advection-vitesse-variable}).
We have done this work in a first numerical experiment, and the results are displayed in Figure \ref{fig-02}.

\smallskip \noindent 
During the first time steps (see the results for $T=01$ and $T=10$), the two results agree with  good precision. 
But we observe that the solution of the advection   equation (\ref{advection-vitesse-variable})
is unsteady, whereas the lattice Boltzmann scheme  rapidly converges  towards a stationary solution.
Then the approximation of the D1Q3 scheme by the first-order partial differential equation
is not sufficient.
We adapt a complementary analysis in the next section. 

\renewcommand{\thefigure}{2}
\begin{figure}    [H]  \centering
\centerline  {\includegraphics[width=.98\textwidth]   {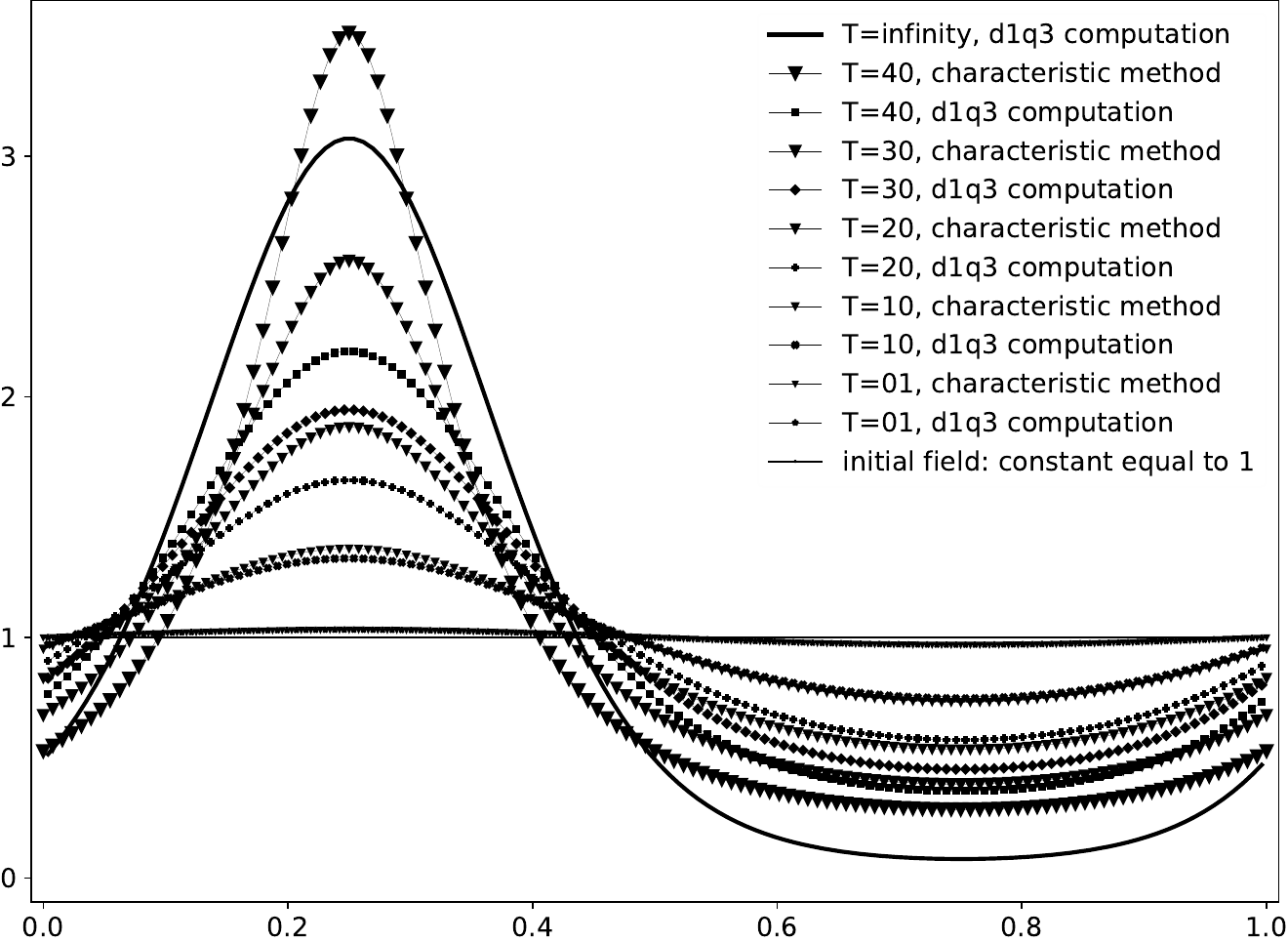}}

\caption{Evolution for an advective velocity (\ref{u-cosinus}) with 
$ \, U = 0.005 \, $ and $ \, N = 128 \,$ mesh points.
We observe that with the cosine advection velocity, the numerical solution has no symmetry. 
As time progresses, we observe that the characteristic  
method becomes less and less precise.}
\label{fig-02} \end{figure}

\bigskip \bigskip    \noindent {\bf \large    4) \quad  ABCD asymptotics in the isotropic linear case} 

\smallskip \noindent 
In this section, we revisit the  ``Berlin algorithm''
introduced in \cite{ADGL14} for the linear analysis of lattice Boltzmann schemes. 
First recall that 
one time step is comprised of two stages:

 \smallskip \noindent 
 {\it (i)} nonlinear relaxation
 \moneqstar
m \, \longmapsto \,  m^* \equiv \begin{pmatrix} W^* \\ Y^* \end  {pmatrix}  \,,\,\, 
 W^* = W   \,,\,\,  Y^* = Y + S \, ( \Phi(W) - Y )
 \monendstar 
 with a  diagonal  relaxation matrix   $ \, S$.
Observe that with our definition, the relaxation matrix~$\, S \, $ only concerns
   the nonconserved microscopic moments. In consequence,
for the D1Q3 scheme, we have $ \, S = {\rm diag} (s\,,\,\, s') $. 
 
\smallskip \noindent 
{\it (ii)}  linear advection
\moneqstar
m^* \, \longmapsto \, f(t+\Delta t):  \,\, f^{*} = M^{-1} \,  m^{*}  \,,\,\, 
f_j(x,\, t+\Delta t) =  f_j^*(x - v_j \, \Delta t , \, t)  .
\monendstar 
The  
operator of advection of moments  $\, \Lambda \, $ is defined 
from the diagonal advection operator $ \,\,  \sum_{\alpha}  v^\alpha\, \partial_\alpha \,\, $ 
according to \cite{Du22}
\moneqstar
\Lambda \equiv
M \,\, {\rm diag} \Big(  \sum_{1  \leq \alpha \leq d}  v^\alpha\, \partial_\alpha \Big) \, M^{-1}  
\monendstar 
with $ \, d \, $ the spatial  dimension. 
This is nothing more  than 
the  advection operator in the basis of moments. 
Following a remark proposed in \cite{Du22}, we have an  
exact exponential expression   of the lattice Boltzmann scheme
\moneqstar
m  (x, t + \Delta t) =  {\rm exp} ( - \Delta t\, \Lambda ) \,\,  m^*(x ,\, t) .
\monendstar 
We can expand this relation up to fourth order to obtain 
\moneqstar
m  (x, t + \Delta t) = \Big[ {\rm I} - \Delta t  \, \Lambda + {{\Delta t^2}\over{2}} \,  \Lambda^2 -  {{\Delta t^3}\over{6}} \,  \Lambda^3
+ {{\Delta t^4}\over{24}} \,  \Lambda^4 + {\rm O}(\Delta t^5) \Big] \, m^*(x ,\, t) .
\monendstar 
The equivalent partial differential equations of the scheme are found  from  the  asymptotic expansion (see {\it e.g.} \cite{DBL23})
\moneqstar
\partial_t = \partial_{t_1}  + \Delta t \, \partial_{t_2} + \Delta t^2 \, \partial_{t_3}
+  \Delta t^3 \, \partial_{t_4} +   {\rm O}(\Delta t^3) . 
\monendstar
Consider now the
``ABCD''  block decomposition (see \cite{DBL23}) of the  momentum-velocity operator that is obtained for the D1Q3 lattice Boltzmann scheme: 
\moneqstar
\Lambda \equiv  \begin{pmatrix} A & B  \\ C  &    D \end  {pmatrix} 
= \left( \begin{array}{c|cc}  0  & \partial_x    & 0  \\  \hline   {2\over3} \, \lambda^2 \, \partial_x   &  0   & \partial_x  \\
  0 &  \lambda^2 \, \partial_x  & 0 \end{array}    \right)  .
\monendstar
%
Asymptotic analysis is carried out to second order. 
It uses only a small set of mathematical expressions: 
\moneq \label{asymptotique-ordre-2} \left\{ \begin{array} {l}
    \partial_{t_1} W +  \Gamma_1 = 0 \\
    \partial_{t_2} W +   \Gamma_2 = 0 \\
    \Gamma_1 =  A \, W + B \, \Phi(W) \\
    Y = \Phi(W) +  \Delta t \,  S^{-1} \, \Psi_1 + {\rm O}(\Delta t^2)  \\
    \Psi_1 =  \dd \Phi(W) .  \Gamma_1  - (C \, W + D \, \Phi(W)) \\
    \Sigma \equiv  S^{-1} - {1\over2} \, {\rm I} \\
    \Gamma_2 =  B \,\,  \Sigma \, \Psi_1  . 
\end{array} \right. \monend 
The application to the Navier Stokes equations can be found  in \cite{dl23}.
For the fourth-order analysis, these relations are enriched in the following way \cite {Du22, DBL23}.
We first extend the 
asymptotic expansion for the microscopic moments 
\moneq \label{moments-microscopiques}
Y = \Phi(W) +   S^{-1} \, \big( \, \Delta t \,  \Psi_1 (W)  + \Delta t^2 \, \Psi_2 (W) + \Delta t^3 \, 
\Psi_3 (W) \, \big) + {\rm O}(\Delta t^4) .
\monend 
We observe that the operators $ \, \Psi_j \, $ are {\it a priori} nonlinear operators of order $ \, j $. 
The partial differential equation for the conserved moments takes the form 
\moneqstar 
  \partial_{t_1} W + \Gamma_1 = 0 \,,\,\, 
  \partial_{t_2} W + \Gamma_2 = 0 \,,\,\, 
  \partial_{t_3} W +  \Gamma_3 = 0 \,,\,\, 
  \partial_{t_4} W +  \Gamma_4 = 0  \, . 
\monendstar 
The differential operators at third order are obtained by nontrivial algebra  \cite {Du22, DBL23}:
\moneq \label{Psi2-Gammma3} \left\{ \begin{array} {l}
  \Psi_2 (W)  =  \Sigma \,  \dd \Psi_1 (W) .   \Gamma_1 (W) + \dd \Phi(W) .  \Gamma_2 (W) - D \, \Sigma \,  \Psi_1 (W)  \\
 \Gamma_3(W) = B \, \Sigma  \, \Psi_2 (W) + {{1}\over{12}}  B_2 \, \Psi_1  (W)
-  {{1}\over{6}} \, B \,  \dd \Psi_1 (W) .  \Gamma_1 (W) 
\end{array} \right. \monend  
and it is also the case for the fourth-order terms:
\moneq \label{Psi3-Gammma4}  \left\{ \begin{array} {l}
  \Psi_3 (W) = \Sigma \, \dd \Psi_1 (W) .  \Gamma_2 (W)  +  \dd \Phi(W) .  \Gamma_3(W) -  D \, \Sigma \, \Psi_2 (W)
  + \Sigma \, \dd \Psi_2 (W) .  \Gamma_1 (W)  \\ \qquad  \qquad
  +{1\over6} \, D \, \dd \Psi_1 (W) .  \Gamma_1 (W)
  - {1\over12} \, D_2 \, \Psi_1 (W)  - {1\over12} \, \dd \, (\dd \Psi_1 (W) .  \Gamma_1 ) .  \Gamma_1 (W) \\
\Gamma_4(W) =  B \, \Sigma \, \Psi_3 (W) + {1\over4} \, B_2 \, \Psi_2 (W)  +  {1\over6} \, B \, D_2 \, \Sigma \, \Psi_1 (W) 
 -   {1\over6} \, A \, B \, \Psi_2 (W) \\ \qquad \qquad 
 -  {1\over6} \, B \, \big( \dd \, (\dd \Phi . \Gamma_1 )  . \Gamma_2 (W)
 -  {1\over6} \, B \, \dd \, (\dd \Phi . \Gamma_2 ) . \Gamma_1 (W) \big) \\ \qquad \qquad
 - {1\over6} \, B \, \Sigma \, \dd \, (\dd \Psi_1 (W) .  \Gamma_1 ) .  \Gamma_1 (W) , 
\end{array} \right. \monend
with
\moneqstar 
\begin{pmatrix} A_2 & B_2  \\ C_2  &    D_2 \end  {pmatrix} \equiv  \begin{pmatrix} A & B  \\ C  &    D \end  {pmatrix}  \, 
\begin{pmatrix} A & B  \\ C  &    D \end  {pmatrix}
= \begin{pmatrix} A^2 + B \, C & A \, B + B \, D \\ C \, A + D \, C & C \, B + D^2 \end  {pmatrix} . 
\monendstar

\bigskip \noindent 
In one spatial  dimension, the previous A, B, C, D differential operators 
factorize  
into algebraic expressions multiplying $x$-derivatives: 
\moneq  \label{Abarre-et-al}
A \equiv {\overline A} \,\, \partial_x \,,\,\,
B \equiv {\overline B} \,\, \partial_x \,,\,\,
C \equiv {\overline C} \,\, \partial_x  \,,\,\,
D \equiv {\overline D} \,\, \partial_x . 
\monend 
\noindent
In this context of
one spatial   
 dimension,  we introduce an  inhomogeneous equilibrium: 
\moneq  \label{Phi-non-uniforme}
\Phi(W) \equiv E(x) \, W \, . 
\monend 
Then we can define a new inhomogeneous  differential operator $ \, \delta \, $ with 
\moneq  \label{operateur-delta}
\delta W  \equiv  \partial_x \big( \Phi(W) \big) =  \partial_x \big( E(x) \, W  \big)  \,. 
\monend
We have  $\, \delta = {{\partial E}\over{\partial x}} \, {\rm I} + E(x) \, \partial_x $.
Then we observe that the commutator $ \, [\partial_x ,\, \delta] \equiv \partial_x \, \delta  - \delta \,  \partial_x \,\, $
is not equal to zero:
$ \,\, [\partial_x \,,\,\, \delta] \, \varphi =
\partial_x \, \partial_x ( E(x) \, \varphi) - \partial_x (E(x) \, \partial_x \varphi) = 
\partial_x \big( (\partial_x E) \, \varphi \big) $. 

\bigskip \noindent
{\bf Proposition 3. Differential operators for linear nonuniform advection}

\noindent
In the previous context of a linear inhomogeneous  scheme, 
the differential operators $ \, \Gamma_1 $, $ \, \Psi_1 $,
$ \, \Gamma_2 $, $ \, \Psi_2 $,  $ \, \Gamma_3 $, $ \, \Psi_3 \, $ and $ \, \Gamma_4 $, 
defined at the relations (\ref{asymptotique-ordre-2})(\ref{Psi2-Gammma3})(\ref{Psi3-Gammma4}),
take the form 
\moneqstar 
\Gamma_j \equiv \alpha_j \, W \,,\,\,  
\Psi_j  \equiv \beta_j \, W 
\monendstar
with the following algebraic relations
\moneq   \label{operateurs-barres}  \left\{ \begin{array} {l}
  \alpha_1 =  {\overline A} \,\, \partial_x + {\overline B} \,\, \delta \\
  \beta_1 = E \, \alpha_1 - ( {\overline C} \,\, \partial_x + {\overline D} \,\, \delta ) \\
  \alpha_2 =  {\overline B} \, \Sigma \, \partial_x \, \beta_1  \\
  \beta_2 = \Sigma \, \beta_1 \, \alpha_1 + E \, \alpha_2 - {\overline D} \,  \Sigma \, \partial_x \, \beta_1 \\ 
   \alpha_3 =  {\overline B} \,  \Sigma \, \partial_x \, \beta_2 
   + {1\over12} \,  {\overline {B_2}} \, \partial_x^2 \,  \beta_1  -  {1\over6} \,  {\overline B} \, \partial_x \, \beta_1 \, \alpha_1 \\
 \beta_3 = \Sigma \, \beta_1 \, \alpha_2 + E \, \alpha_3 -  {\overline D} \,  \Sigma \, \partial_x \, \beta_2
 +  \Sigma \, \beta_2 \, \alpha_1 + {1\over6} \, {\overline D} \,   \partial_x \, \beta_1 \, \alpha_1
 -{1\over12} \,  \beta_1 \, \alpha_1^2  -{1\over12} \,  {\overline D_2} \,   \partial_x^2 \, \beta_1 \\  
 \alpha_4 =  {\overline B} \,  \Sigma \, \partial_x \, \beta_3 +  {1\over4} \,  {\overline B_2} \, \partial_x^2 \, \beta_2
+  {1\over6} \, {\overline B} \, {\overline D_2} \,  \Sigma \,\partial_x^3 \,  \beta_1
-  {1\over6} \, {\overline A} \, {\overline B} \,\partial_x^2 \,  \beta_2 \\ \qquad 
-  {1\over6} \,{\overline B} \, \delta \, \alpha_1 \, \alpha_2
-  {1\over6} \,{\overline B} \, \delta \, \alpha_2 \, \alpha_1  -  {1\over6} \,{\overline B} \,  \Sigma \, \partial_x \, \beta_1 \, \alpha_1^2 \, . 
\end{array} \right. \monend
The proof of this proposition is a tedious algebraic calculation. It is presented in Annex A. 

\smallskip \smallskip  \noindent 
We consider in this contribution  the case of the D1Q3 scheme with one conservation law with 
a cosine velocity field  $ \, u(x) \equiv  \lambda \, U \, \cos(k \, x) \, $ introduced in (\ref{u-cosinus}). 
Then  the differential operator $ \, \delta \, $ proposed in~(\ref{operateur-delta}) 
takes the form
\moneqstar
\delta \varphi = \partial_x \big( E(x) \, \varphi \big)
= \partial_x \begin{pmatrix} \lambda \, U \, \cos(k \, x) \varphi \\ \lambda^2 \, \alpha \, \varphi \end{pmatrix} 
= \partial_x \begin{pmatrix} u \, \varphi \\ \lambda^2 \, \alpha \, \varphi \end{pmatrix} . 
\monendstar 
With the notation
\moneq   \label{partial_u}
\partial_u \varphi \equiv U \, \partial_x \big(\cos(k \, x) \varphi  \big)
\monend
if the velocity field is a cosine  
({\it c.f.} (\ref{u-cosinus})) and
\moneqstar 
\partial_u \varphi \equiv U \, \partial_x \varphi  
\monendstar
when the velocity field is constant, 
we have simply 
\moneq   \label{delta-plus-simple}
\delta =  \begin{pmatrix} \lambda  \, \partial_u \\ \lambda^2 \, \alpha \, \partial_x \end{pmatrix} .
\monend
We observe  that 
the differential operators $\,  \partial_x \, $ and $\,  \partial_u \, $ do not commute:
\moneqstar
 [ \partial_x \,,\,\,  \partial_u ] \, \varphi = \partial_x \big( - k \, U \, \sin (k\,x) \,\, \varphi \big) . 
\monendstar 
%

\bigskip \noindent
{\bf Proposition 4. D1Q3 differential operators for linear nonuniform advection}

\noindent
The linear differential operators explicited  in  (\ref{operateurs-barres})
admit the following expressions in terms of the operators $ \, \partial_x \, $ and $ \, \partial_u $.
We have for the two first orders 
\moneq   \label{alpha1-beta1} 
\alpha_1 = \lambda \, \partial_u  \,,\,\, 
 \beta_1 = \lambda^2 \,\begin{pmatrix}   {{u}\over{\lambda}} \, \partial_u  - {{\alpha+2}\over{3}} \, \partial_x   \\
 \lambda \, (\alpha - 1) \,  \partial_u \end{pmatrix} \monend
%
\moneq   \label{alpha2-beta2} \left\{ \begin{array} {l}
  \alpha_2 = \lambda^2 \, \sigma \,  \Big(  \partial_u^2 -  {{\alpha+2}\over{3}} \, \partial_x^2 \Big) \\ 
 \beta_2 = \lambda^3 \, \begin{pmatrix} 2 \, \sigma \, {{u}\over{\lambda}} \,  \partial_u^2 
  - \big( {{\alpha+2}\over{3}}\, \sigma + {{\alpha-1}\over{3}}\, \sigma' \big)\, \partial_x \, \partial_u 
  - {{\alpha+2}\over{3}}\, \sigma \, {{u}\over{\lambda}} \, \partial_x^2  \\
  \lambda \,(\alpha - 1) \, \big( (\sigma + \sigma')\, \partial_u^2 
  - {{\alpha+2}\over{3}}\, \sigma \, \partial_x^2 \big) \end{pmatrix} ,
\end{array} \right. \monend
with the H\'enon coefficients \cite{He87} 
$\, \sigma \, $ and $ \, \sigma' \, $ defined according to
\moneqstar
\sigma = {1\over{s}} - {1\over2} \,,\,\, \sigma' = {1\over{s'}} - {1\over2} .
\monendstar 
At third order, the formulae are more complicated. We have 
\moneq   \label{alpha3}
\alpha_3 = \lambda^3  \Big[ \Big( 2 \, \sigma^2 - {1\over6} \Big) \, \partial_u^3 
  + \Big(  {{\alpha+2}\over{3}}\, \big( {1\over6}-\sigma^2 \big) + {{\alpha-1}\over{3}}\,  \big( {1\over12}- \sigma\, \sigma' \big)
  \Big) \,  \partial_x^2 \,\partial_u    - {{\alpha+2}\over{3}}\,\sigma^2 \, \partial_u \, \partial_x^2   \Big]
\monend 
and
\moneqstar 
\beta_3 \equiv \begin{pmatrix} \lambda^4 \, \beta_{3J} \\  \lambda^5 \, \beta_{3e} \end{pmatrix} ,
\monendstar
with
\moneq   \label{beta3} \left\{ \begin{array} {l}
\beta_{3J} = {{\alpha+2}\over9}\, \big[ -(1-\alpha)\, \sigma \, \sigma'
  +  (\big (\alpha +2)\, \sigma^2  + {1\over4} \big) \big] \, \partial_x^3 \\ \qquad 
+ U \,  \big[  -2 \,  {{\alpha+2}\over3}\, \sigma^2 + {{1-\alpha}\over3} \,\sigma \, \sigma' +   {{1+\alpha}\over12} \big] 
\, \partial_x^2  \, \partial_u - 2\, U \, {{\alpha+2}\over3} \, \sigma^2 \, \partial_u \, \partial_x^2  \\ \qquad 
 +  \big[ -2\,  {{\alpha+2}\over3}\, \sigma^2
  +  {{1-\alpha}\over3} \,(2 \, \sigma \, \sigma' + \sigma'^2 -{1\over4}) \big] \, \partial_x \, \partial_u^2 
+ \big( 5 \, \sigma^2  -{1\over4} \big) \, U \,  \partial_u^3  \\ 
 \beta_{3e} =  {{1-\alpha}\over3} \,  \big[  (\alpha+2) \, \sigma^2 + (1 + 2\, \alpha) \, \sigma \, \sigma'
  -  {{1+\alpha}\over4} \big] \, \partial_x^2 \, \partial_u 
+ (1-\alpha) \,  {{\alpha+2}\over3}\, \sigma \, ( \sigma +  \sigma' ) \, \partial_u \, \partial_x^2  \\  \qquad 
- (1-\alpha) \, \big( 2 \,  \sigma^2 + 2 \, \sigma \, \sigma' + \sigma'^2  -{1\over4}  \big) \, \partial_u^3 .
\end{array} \right. \monend
At fourth order, we have 
\moneq   \label{alpha4} \left\{ \begin{array} {l}
\alpha_4 = \lambda^4 \, \Big[ \,
  {{\alpha+2}\over9} \, \big( (\alpha+2) \, \sigma^3 - (1-\alpha)\, \sigma^2 \, \sigma'
  - {{\alpha}\over4} \, \sigma \big) \,  \partial_x^4 \\ \qquad 
  + \big[ -2\, {{\alpha+2}\over3} \, \sigma^3
+  {{1-\alpha}\over3} \, (2 \, \sigma^2 \, \sigma' + \sigma \, \sigma'^2 - {1\over4} \, \sigma' )
+ {{1 + 2\, \alpha}\over9} \,   \sigma \big] \,  \partial_x^2 \, \partial_u^2  \\  \qquad 
  + \big[ -2\, {{\alpha+2}\over3} \, \sigma^3 +   {{1-\alpha}\over3} \, \sigma^2 \, \sigma'
    + {{7 + 5 \, \alpha}\over36} \, \sigma \big] \,  \partial_u \, \partial_x^2 \, \partial_u  \\  \qquad 
  + {{\alpha+2}\over3} \, \sigma \, \big( -2 \, \sigma + {1\over6} \big) \, \partial_u^2 \, \partial_x^2 
  + \sigma \, \big( 5 \, \sigma^2 - {3\over4} \big) \,  \partial_u^4  \,  \Big] \, . 
\end{array} \right. \monend

The proof of Proposition 4 is detailed in  Annex B. 

\bigskip \noindent
When the velocity field has a constant value $ \, u(x) \equiv \lambda \, U $, 
super-convergence can be obtained with an appropriate choice of relaxation coefficients,
called ``magic'' in \cite{HG09}. Because magic is not science,  we prefer the
denomination of ``quartic parameters''  \cite{DL11} to achieve fourth-order accuracy,
or  ``cubic parameter'' in the present case to obtain a third-order precision. 

\smallskip \noindent
When the advection velocity field has a constant value, we have $ \, \partial_u \equiv U \, \partial_x \, $
and the coefficient~$ \, \alpha_3 \, $ initially given according to (\ref{alpha3}) takes now the value 
\moneq   \label{alpha3bis}
\alpha_3 = {{\lambda^3 \, U}\over12} \, \Big[ -2 \, (1-12\, \sigma^2)\, U^2 + 4 \,(1-\alpha) \, \sigma \, \sigma' + 1 + \alpha
 -8 \, (2+\alpha) \, \sigma^2 \Big] \, \partial_x^3  \, . 
\monend 
Then for a fixed set of values for  $ \, U $, $ \, \alpha \, $ and $ \, \sigma $,  the cubic  parameter
$ \, \sigma_c' \, $  is defined  by forcing to zero the value of $ \, \alpha_3 \, $ in the relation (\ref{alpha3bis}): 
\moneq   \label{sigma-prime-cubique}
\sigma_c' = {{2 \, (1-12\, \sigma^2)\, U^2 + 8 \, (2+\alpha) \, \sigma^2 - (1+\alpha)}\over
  {4 \,(1-\alpha) \, \sigma}} \, . 
\monend 

\smallskip \smallskip  \noindent
In the following, we first experiment with the D1Q3 lattice Boltzmann scheme with constant velocity, 
possibly with cubic parameters. Then we consider a cosine advection velocity.
We detail in the next section the Fourier methodology developed to solve
the various equivalent partial differential equations with very good precision.

\bigskip \bigskip 
\noindent {\bf \large    5) \quad  Interlaced Fourier series} 

\smallskip \noindent
We  compare the numerical simulation done with the D1Q3 lattice Boltzmann scheme
with the solution of the  equivalent partial differential equations up to fourth-order accuracy.
This hierarchy of equations can be  written 
\moneq   \label{equivalent-edp}
{{\partial \rho}\over{\partial t}}  +  \sum_{j=1}^\ell \Delta t^{j-1} \, \alpha_j \,  \rho  = 0 .
\monend
They  are of order $ \, \ell \, $ for $ \, 1 \leq \ell \leq 4 $.
We recall that we have the following structure
\moneqstar \left\{ \begin{array} {l}
\Delta t^0 \, \alpha_1 = \lambda \, \partial_u   \\ 
 \Delta t^1 \,  \alpha_2 = -\mu \, \partial_x^2 + \mu_u \, \partial_u^2  \\ 
  \Delta t^2 \,  \alpha_3 = \xi_u \, \partial_u^3 + \xi_{xu}  \,\partial_x^2 \, \partial_u  + \xi_{ux}   \, \partial_u \,\partial_x^2   \\
  \Delta t^3 \,  \alpha_4 =  \zeta_{u4} \, \partial_u^4 + \zeta_{xxuu}  \, \partial_x^2 \, \partial_u^2
  +  \zeta_{uxxu}  \, \partial_u \, \partial_x^2 \, \partial_u  + \zeta_{uuxx}  \, \partial_u^2 \, \partial_x^2 +  \zeta_{x4} \, \partial_x^4 .
\end{array} \right. \monendstar 
The coefficients $ \, \mu $, $ \, \mu_u $, $ \, \xi_u $,  $ \, \xi_{xu} $, $ \, \xi_{ux} $,
$ \, \zeta_{u4} $, $ \, \zeta_{xxuu} $, $ \,  \zeta_{uxxu} $, $ \, \zeta_{uuxx} \,$ and $ \,  \zeta_{x4}  \, $ 
are easy to explicate from the relations (\ref{alpha1-beta1})(\ref{alpha2-beta2})(\ref{alpha3})(\ref{alpha4}):
%
\moneqstar \left\{ \begin{array} {l}
  \mu = {{\alpha+2}\over{3}} \, \lambda^2 \, \sigma \,,\,\, \mu_u = \lambda^2 \, \sigma \,, \\
  \xi_u =  \lambda^3 \,  \big( 2 \, \sigma^2 - {1\over6} \big)  \,,\,\,
  \xi_{xu}  =  \lambda^3 \, \Big(  {{\alpha+2}\over{3}}\, \big( {1\over6}-\sigma^2 \big)
  + {{\alpha-1}\over{3}}\,  \big( {1\over12}- \sigma\, \sigma' \big)   \Big)  \,,\,\,
  \xi_{ux} = -  \lambda^3 \,  {{\alpha+2}\over{3}}\,\sigma^2  \\
\zeta_{u4} =  \lambda^4 \,  \sigma \, \big( 5 \, \sigma^2 - {3\over4} \big)  \,,\,\,
 \zeta_{xxuu} =  \lambda^4 \, \big[ -2\, {{\alpha+2}\over3} \, \sigma^3
+  {{1-\alpha}\over3} \, (2 \, \sigma^2 \, \sigma' + \sigma \, \sigma'^2 - {1\over4} \, \sigma' )
+ {{1 + 2\, \alpha}\over9} \,   \sigma \big] \\
 \zeta_{uxxu} =  \lambda^4 \, \big[ -2\, {{\alpha+2}\over3} \, \sigma^3 +   {{1-\alpha}\over3} \, \sigma^2 \, \sigma'
   + {{7 + 5 \, \alpha}\over36} \, \sigma \big]   \,,\,\,
 \zeta_{uuxx} =  \lambda^4 \, {{\alpha+2}\over3} \, \sigma \, \big( -2 \, \sigma + {1\over6} \big) \\ 
\zeta_{x4} = \lambda^4 \, \big[ \,   {{\alpha+2}\over9} \, \big( (\alpha+2) \, \sigma^3 - (1-\alpha)\, \sigma^2 \, \sigma'
  - {{\alpha}\over4} \, \sigma \big) \big] . \end{array} \right. \monendstar 
We use a spectral method to capture an approximation of a partial differential equation of the family (\ref{equivalent-edp}).
In the case of an advective field given in (\ref{u-cosinus}), we introduce the  two discrete spaces
$ \, S_i \, $ and $ \, S_p \, $ defined as follows.  
The space of odd sine  and even cosine  is called  $ \, S_i $:
\moneqstar
S_i \ni \rho = 
\sum_{j \geq 0} a_{2 j + 1} \, \sin \big( (2 j + 1)\,k \,x \big)
+ \sum_{j \geq 0} a_{2 j + 2} \, \cos \big( (2 j + 2)\,k \,x \big)
\monendstar
and the space of even   sine and odd cosine  is denoted by  $ \, S_p $:
\moneqstar
S_p \ni \rho = \sum_{j \geq 0} b_{2 j + 1} \, \cos \big( (2 j + 1)\,k \,x \big)  + 
\sum_{j \geq 0} b_{2 j + 2} \, \sin \big( (2 j + 2)\,k \,x \big) .
\monendstar
The derivation operator breaks down into two parts:
\moneqstar \left\{ \begin{array} {l}
  \partial_x^{ip} \, : \, S_i \longrightarrow S_p  \\
  \partial_x^{pi} \, : \, S_p \longrightarrow S_i \, .
\end{array} \right. \monendstar
Relatively to the basis
$ \, \big( \sin k\,x ,\, \cos \, 2 \, k\,x  ,\, \sin \, 3 \, k\,x ,\,  \cos \, 4 \, k\,x  ,\, \cdots \big) \, $ of  $ \, S_i \, $
and  to the basis
$ \, \big( \cos k\,x ,\, \sin \, 2 \, k\,x  ,\, \cos \, 3 \, k\,x ,\,  \sin \, 4 \, k\,x  ,\, \cdots \big) \, $ of  $ \, S_p  $,  
the operators  $ \, \partial_x^{ip} \, $ and $ \, \partial_x^{pi} \, $ 
can be represented by the following matrices 
\moneqstar \left\{ \begin{array} {l}
  \Delta_x^{ip} = {\rm diag} \big( k ,\, -2\, k ,\, 3 \, k ,\, -4\, k ,\, \cdots \big)  \\
  \Delta_x^{pi} = {\rm diag} \big( -k ,\, 2\, k ,\, -3 \, k ,\, 4\, k ,\, \cdots \big) = -  \Delta_x^{ip} . 
\end{array} \right. \monendstar

\smallskip \noindent 
The second order operator $ \, \partial_x^2 = \partial_x^{pi} \,_\circ \,  \partial_x^{ip} \, $
operates inside the space $ \, S_i \, $ and is represented by the matrix
\moneqstar
  \Delta_x^{pi}  \,\,    \Delta_x^{ip} = -  {\rm diag} \big( k^2 ,\, 4\, k^2 ,\, 9 \, k^2 ,\, 16 \, k^2 ,\, \cdots \big) .
\monendstar

\smallskip \noindent 
We introduce also  the operator $ \, m_u \, $ of mutiplication by $ \, u \equiv U \, \cos(k\, x) $.
It operates from  $\, S_i \, $ and takes its values in $ \, S_p $.
Then 
$ \, \partial_u = \partial_x \, _\circ\, m_u \, $ operates inside the space  $\, S_i   $.
More precisely, we have, without forgetting the  constant component $ \, a_0 $: 

\smallskip 
$  m_u \, \rho = U \, \cos(k\, x) \, \Big[ a_0 + \sum_{j \geq 0} \, a_{2 j + 1} \, \sin \big( (2 j + 1)\,k \,x \big)
  + \sum_{j \geq 0} \, a_{2 j + 2} \, \cos \big( (2 j + 2)\,k \,x \big) \Big] $

\smallskip \qquad 
$ \, \, = \big(a_0 + {1\over2}\, a_2 \big) \, U \,  \cos \big( k \,x \big)
+ {{U}\over2}\, \sum_{j \geq 0} \big( a_{2 j + 1} +  a_{2 j + 3} \big) \,  \sin \big( (2  j + 2)\,k \,x \big) $

\smallskip \qquad \qquad 
$ + {{U}\over2}\, \sum_{j \geq 1}  \big( a_{2 j} +  a_{2 j +2} \big) \,\cos \big( (2  j + 1)\,k \,x \big) .  $


\smallskip \noindent
The matrix
\moneqstar
  M_u \equiv {{U}\over2}\,  \begin {pmatrix} 0 & 1 & 0 &0 \\ 1 & 0 & 1 & 0 \\ 0 & 1 & 0 & \ddots \\ 0 & 0 & \ddots & \ddots    \end{pmatrix}
\monendstar 

\smallskip \noindent
is a  natural implementation of the operator $ \, m_u \, $ of multiplication  by the velocity $\, u \, $
for $\, \rho \in  S_i $. 
Then the differential operator
$ \, \partial_u \equiv  \partial_x \, m_u \, $ in the space $ \, S_i \, $ after truncation
is represented by the matrix $ \,  \Delta_x^{pi} \, M_u $. 


\smallskip \noindent
We have used two discretizations with 30 or 60 Fourier modes.
Observe that when $ \, U = 0.05 $,
the results are correct with 30 Fourier modes for 64 and 128 mesh points.
But with 256 and 512 mesh points, oscillations appear in the numerical results.
This is the sign of a under-resolved simulation.
We have changed the number of Fourier modes and used 60 modes for  256 and 512 mesh points.
We tested  the representation of the solution of the D1Q3 scheme with a Fourier series.
We have observed a residual in $ \, \ell^\infty \, $ norm of
$ \,  1.79 \times 10^{-14} \, $ and $ \,  6.06 \times 10^{-11} $.
This precision is sufficient for our simulations. 

\smallskip \noindent
After these algebraic operations, the partial differential equation
(\ref{equivalent-edp}) can be seen as  an infinite system of ordinary differential equations
\moneq   \label{edo-infini}
{{\partial \rho}\over{\partial t}}  +  A \, \rho  = 0 
\monend 
with an operator $ \, A \, $ given at fourth order by the relation
\moneq   \label{operateur-A} \left\{ \begin{array} {l}
  A = \lambda \, \partial_x \, m_u   -\mu \, \partial_x^2 + \mu_u \, \partial_x \, m_u \, \partial_x \, m_u \\ \qquad 
  + \big( \xi_u \, \partial_x \, m_u \, \partial_x \, m_u \, \partial_x \, m_u
  + \xi_{xu}  \,\partial_x^3 \,  m_u   + \xi_{ux}   \, \partial_x \, m_u \, \partial_x^2 \big) \\ \qquad 
  +  \big( \zeta_{u4} \, \partial_x \, m_u \, \partial_x \, m_u \, \partial_x \, m_u \, \partial_x \, m_u 
  + \zeta_{xxuu}  \, \partial_x^3 \,  m_u \, \partial_x \, m_u
  +  \zeta_{uxxu}  \, \partial_x \, m_u \, \partial_x^3 \, \, m_u \\  \qquad 
+ \zeta_{uuxx}  \,  \partial_x \, m_u \, \partial_x \, m_u  \, \partial_x^2 +  \zeta_{x4} \, \partial_x^4 \big)   \, .
\end{array} \right. \monend
We observe that the matrix $ \, A \, $ is constant. Then after discretization with $ \, N \, $ modes,
it becomes a constant matrix $ \, \widetilde{A}_N $.
The system (\ref{edo-infini}) is replaced by a system of a finite number of  ordinary differential equations
\moneq   \label{edo-fini}
{{\partial \rho}\over{\partial t}}  +  \widetilde{A}_N \, \rho  = 0 \, . 
\monend 
Due to the fact that the matrix $ \, \widetilde{A}_N \, $ is fixed, the solution of  (\ref{edo-fini}) is approached
in this contribution
for small values of time ($t$ is replaced by the time step $ \, \Delta t \, $ for the numerical integration)
by a Taylor expansion at order 5 from the initial condition $ \, \rho_0  $, 
an expansion consistent  with a  
spatial precision of order at most 4:
\moneqstar
\rho(t) = \exp(-t \, \widetilde{A}_N) \, \rho_0 \simeq \Big[ {\rm I} - t \, \widetilde{A}_N + {{t^2}\over{2}} \, \widetilde{A}_N^2
  - {{t^3}\over{6}} \, \widetilde{A}_N^3  + {{t^4}\over{24}} \, \widetilde{A}_N^4 - {{t^5}\over{120}} \, \widetilde{A}_N^5  \Big]  \, \rho_0 \, . 
\monendstar 

\smallskip \noindent
In an initial  series of numerical experiments, we have put in evidence approximations of the stationary 
solution of a lattice Boltzmann scheme forced with a cosine velocity field.

\bigskip \bigskip    \noindent {\bf \large    6) \quad  Unsteady  evolution} 

\smallskip \noindent 
We  now compare the D1Q3 lattice Boltzmann scheme up to time $ \, T = 1 \, $ with the Fourier approximations of the
various equivalent partial differential equations at various orders 
\moneqstar
{{\partial \rho}\over{\partial t}}  +  A_j \, \rho  = 0  \,, \,\,  1 \leq j \leq 4 .
\monendstar 
Observe that $\,  A_j \,$ is a continuous partial differential operator whereas $ \, \widetilde{A}_N \, $
in the previous section  is a discretization with $ \, N \, $ degrees of freedom. We use a different notation 
to avoid ambiguity between  the order of approximation and the order of discretization.
We have at order 1:
\moneqstar
A_1 = \lambda \, \partial_x \, m_u \,, 
\monendstar 
at order 2: 
\moneqstar
A_2 = A_1    -\mu \, \partial_x^2 + \mu_u \, \partial_x \, m_u \, \partial_x \, m_u \,, 
\monendstar 
at order 3:
\moneqstar
A_3 = A_2 +   \big( \xi_u \, \partial_x \, m_u \, \partial_x \, m_u \, \partial_x \, m_u
+ \xi_{xu}  \,\partial_x^3 \,  m_u   + \xi_{ux}   \, \partial_x \, m_u \, \partial_x^2 \big)  \,, 
\monendstar 
and at order 4: 
\moneqstar \left\{ \begin{array} {l}
A_4 = A_3 +  \big( \zeta_{u4} \, \partial_x \, m_u \, \partial_x \, m_u \, \partial_x \, m_u \, \partial_x \, m_u 
+ \zeta_{xxuu}  \, \partial_x^3 \,  m_u \, \partial_x \, m_u \\ 
\qquad \qquad  +  \zeta_{uxxu}  \, \partial_x \, m_u \, \partial_x^3 \, \, m_u
+ \zeta_{uuxx}  \,  \partial_x \, m_u \, \partial_x \, m_u  \, \partial_x^2 +  \zeta_{x4} \, \partial_x^4 \big) \, . 
\end{array} \right. \monendstar
We have chosen the following parameters
\moneqstar
U = 0.005 \,,\,\,  \alpha = -1 \,,\,\, \sigma = 0.01 \,\,\, [s = 1.960784313725] \,,\,\, s' = 1.2  
\monendstar 
with $ \, \sigma \equiv {1\over{s}} - {1\over2} $.
  The detailed results are presented in a preliminary edition of this work \cite{BDL24}.

\renewcommand{\thefigure}{3}
\begin{figure}    [H]  \centering
  \centerline  {\includegraphics[width=.58\textwidth]   {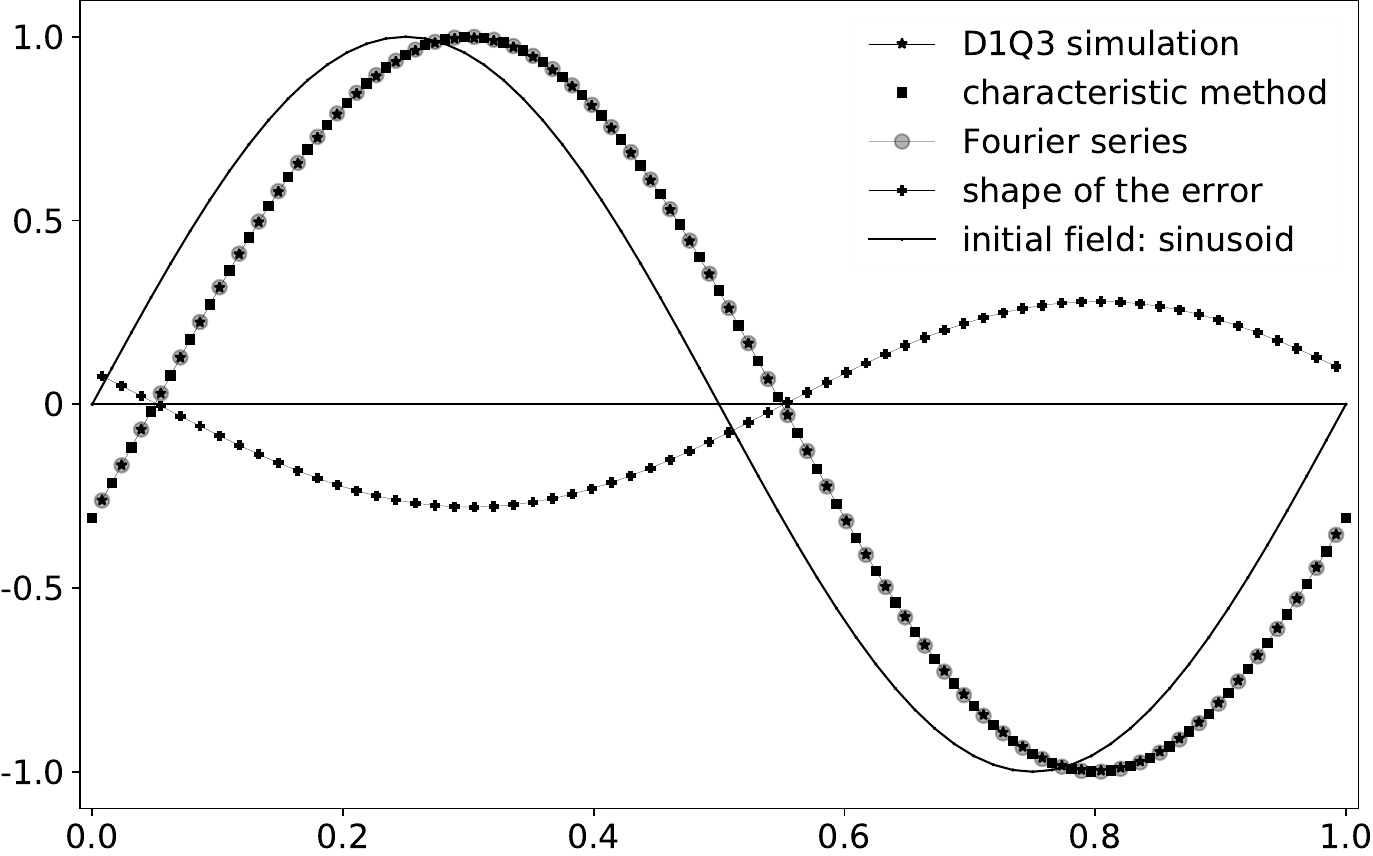}}
\caption{Unsteady evolution, constant advection field [$ U = 0.05  $], 64 mesh points, sinusoidal initial condition.}
\label{fig-06a} \end{figure} 

\renewcommand{\thefigure}{4}
\begin{figure}    [H]  \centering
  \centerline  {\includegraphics[width=.58\textwidth]   {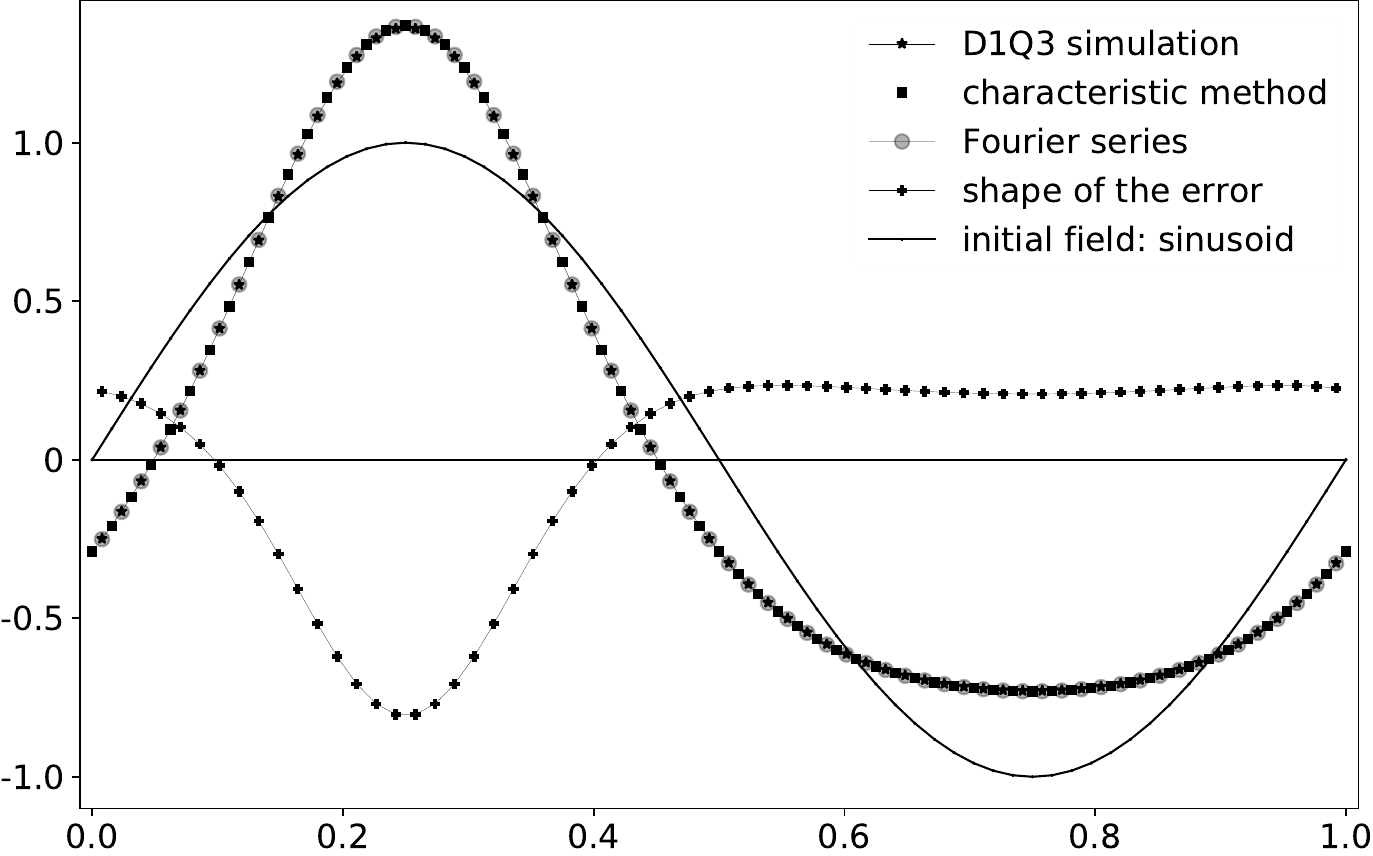}}
\caption{Unsteady evolution, sinusoidal advection field [$ U = 0.05  $], 64 mesh points, sinusoidal initial condition.}
\label{fig-07a} \end{figure} 

%
%
In Figure \ref{fig-06a}, 
we consider a constant velocity field with a sinusoidal initial condition.  
In Figure \ref{fig-07a},  
a cosine  velocity field with a sinusoidal initial condition is studied. 
These figures show that the approximation of the lattice Boltzmann scheme with the equivalent partial
differential equations is globally correct.
Then we refine the mesh up to $ \, N = 1024 \, $ points. 
The results are  presented in Figure  \ref{fig-09}. They are not completely satisfactory.

\renewcommand{\thefigure}{5}
\begin{figure}    [h]  \centering
\centerline  {\includegraphics[width=.58\textwidth]   {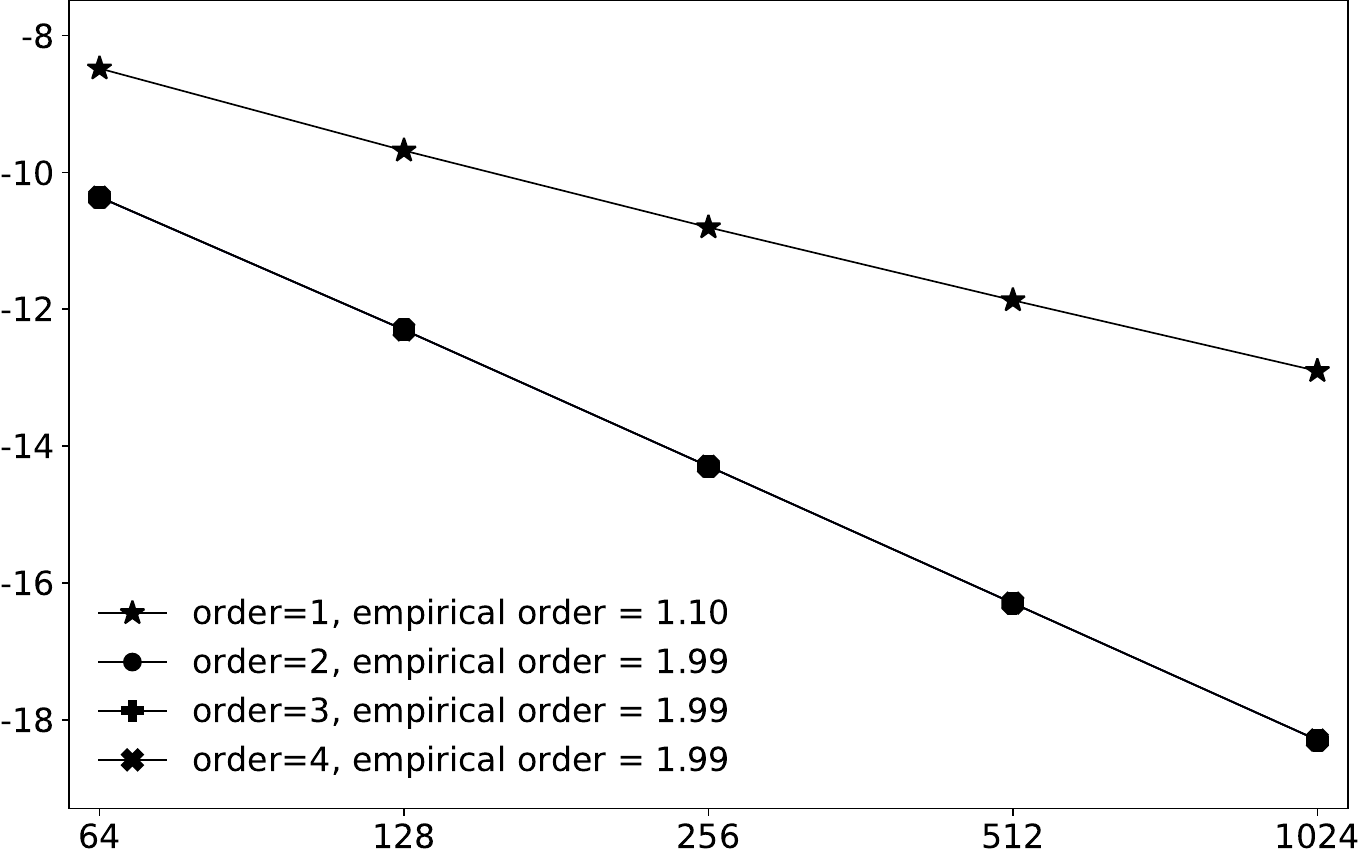}}
\caption{Errors measured with the maximum norm
  between the D1Q3 lattice Boltzmann  scheme
and various equivalent  partial differential  equations. Unsteady experiment with 
constant velocity field $\, U = 0.05$,  finite-time evolution with final time $ \, T = 1 $,
and initialization with a sinusoidal  wave.
The $x$-axis represents the number of mesh points with a logarithmic scale and
the $y$-axis is graduated with the base-2   logarithm  of the error.
The microscopic moments are  initialized with the equilibrium values.}
\label{fig-09} \end{figure} 

\smallskip \noindent
In order to overcome the moderate speed of convergence
for an unsteady evolution, we focus in the next section on the way the lattice Boltzmann scheme
is initialized.

\bigskip \bigskip    \noindent {\bf \large    7) \quad  Initialization of microscopic moments} 

\smallskip \noindent
In the previous section, we have taken the non-conserved moments at time $ \, t = 0 \, $
equal to the value at equilibrium:
\moneq   \label{initialisation-ordre-0}
Y_0(t=0) = \Phi(\rho_0) .
\monend 
We recall  the asymptotic expansion of nonconserved moments for a lattice Boltzmann scheme
through the general expression (\ref{moments-microscopiques}): 
\moneqstar
 Y = \Phi(W) +   S^{-1} \, \big( \, \Delta t \,  \Psi_1 (W)  + \Delta t^2 \,  \Psi_2 (W) \, \big) + {\rm O}(\Delta t^3) .
\monendstar
For the advective D1Q3 scheme with a cosine advection field, 
\moneqstar
Y  =  \begin{pmatrix} j \\ e \end{pmatrix} \,,\,\, 
\Phi(\rho) = \begin{pmatrix} \lambda \, U \, \cos(k\, x) \, \rho \\ \lambda^2 \, \alpha \, \rho \end{pmatrix}  \,,\,\, 
\Psi_1 = \beta_1 \, \rho  \,,\,\,
\beta_1 = \lambda^2 \,\begin{pmatrix}   u \, \partial_u - {{\alpha+2}\over{3}} \, \partial_x  \\
\lambda \, (\alpha - 1) \,  \partial_u \end{pmatrix} 
\monendstar
and 
\moneqstar
\Psi_2 = \beta_2 \, \rho  \,,\,\,
\beta_2 = \lambda^3 \, \begin{pmatrix} 2 \, \sigma \, {{u}\over{\lambda}} \,  \partial_u^2
- \big( {{\alpha+2}\over{3}}\, \sigma + {{\alpha-1}\over{3}}\, \sigma' \big)\, \partial_x \, \partial_u
- {{\alpha+2}\over{3}}\, \sigma \,  {{u}\over{\lambda}} \,\partial_x^2  \\
\lambda \,(\alpha - 1) \, \big( (\sigma + \sigma')\, \partial_u^2 - {{\alpha+2}\over{3}}\, \sigma \, \partial_x^2 \big) 
\end{pmatrix} .
\monendstar
For initialization  at order 0, the relation (\ref{initialisation-ordre-0})
is simply applied. The initialization suggested by Mei, Luo, Lallemand and  d'Humi\`eres \cite{MLLH06}  
at order 1 is: 
\moneq   \label{initialisation-ordre-1}
Y_1(t=0) = \Phi(\rho_0) +  \Delta t \, S^{-1} \, \beta_1 \, \rho_0 .
\monend 
In the following, we also consider  a second-order initialization:
\moneq   \label{initialisation-ordre-2}
Y_2(t=0) = \Phi(\rho_0) +  S^{-1} \, \big( \Delta t \, \beta_1 \, \rho_0   + \Delta t^2 \, \beta_2 \, \rho_0  \big) .
\monend 

\smallskip \noindent 
We remark that this framework can certainly be  revisited with the  new 
version  of lattice Boltzmann schemes through  multistep finite difference schemes, 
as proposed by Bellotti, Graille and Massot  \cite{BGM22}.
The results of our simulations are presented in the next section.

\bigskip \bigskip    \noindent {\bf \large    8) \quad  Unsteady fields for a constant or variable advective velocity} 

\smallskip \noindent 
We first study the uniform advection case. 
Then we specify the case of cubic  parameters.
Then we look to nonuniform cosine advection.
In all cases, the choice of the initialization scheme has a great influence on the final precision. 
Observe also that only one  mode is needed for the Fourier approximation when
the advecting velocity  is constant.

\renewcommand{\thetable}{1}
\begin{table}[H]
\smallskip
\centerline {\begin{tabular}{|c|c|c|c|c|}    \hline 
mesh points $ \backslash  $ equation order &  1  &  2 &  3 &  4 \\   \hline
initialization order  & 0 & 0 & 0 & 0  \\  \hline
  64 & 2.798 $10^{-3}$  & 7.606 $10^{-4}$  & 7.604 $10^{-4}$  & 7.596 $10^{-4}$   \\  \hline
 128 & 1.218 $10^{-3}$  & 1.983 $10^{-4}$  & 1.983 $10^{-4}$  & 1.982 $10^{-4}$   \\  \hline
 256 & 5.598 $10^{-4}$  & 4.979 $10^{-5}$  & 4.979 $10^{-5}$  & 4.978 $10^{-5}$   \\  \hline
 512 & 2.675 $10^{-4}$  & 1.245 $10^{-5}$  & 1.245 $10^{-5}$  & 1.245 $10^{-5}$   \\  \hline
1024 & 1.307 $10^{-4}$  & 3.113 $10^{-6}$  & 3.112 $10^{-6}$  & 3.112 $10^{-6}$   \\  \hline
convergence order & 1.10 & 1.99 & 1.99 & 1.99  \\  \hline
\end{tabular}}
\caption{Errors measured with the maximum norm
between the D1Q3 lattice Boltzmann  scheme [with parameters 
$ \, \alpha = -1 $,  $ \, \sigma \equiv {{1}\over{s}} - {1\over2} = 0.01$, $\, s' = 1.2 $]  
and various equivalent  partial differential  equations for an unsteady experiment:
constant velocity field $\, U = 0.05$,  finite time evolution with final time $ \, T = 1 $,
and initialization with a sinusoidal wave.
The error remains second-order accurate even if we use the third-order or the fourth-order
equivalent  equation for the approximation of the lattice Boltzmann scheme.
Figure \ref{fig-09} is an other representation of these results.}
\label{tab-04} \end{table}

\smallskip \noindent
With the first-order initialization  (\ref{initialisation-ordre-1}), the results are presented in
Figure \ref{fig-10} and Table  \ref{tab-05}.
They become consistent for the three first levels of approximation, but
there is no convergence at fourth-order accuracy. 

\smallskip \noindent 
With the second-order initialization  (\ref{initialisation-ordre-2}), the results are displayed
in Figure \ref{fig-11} and Table~\ref{tab-06}.
The experimental order of approximation is now coherent up to fourth order.
In Table  \ref{tab-07}, economical  initialization orders are used to present an optimal convergence accuracy.

\renewcommand{\thefigure}{6}
\begin{figure}    [h]  \centering
\centerline  {\includegraphics[width=.58\textwidth]   {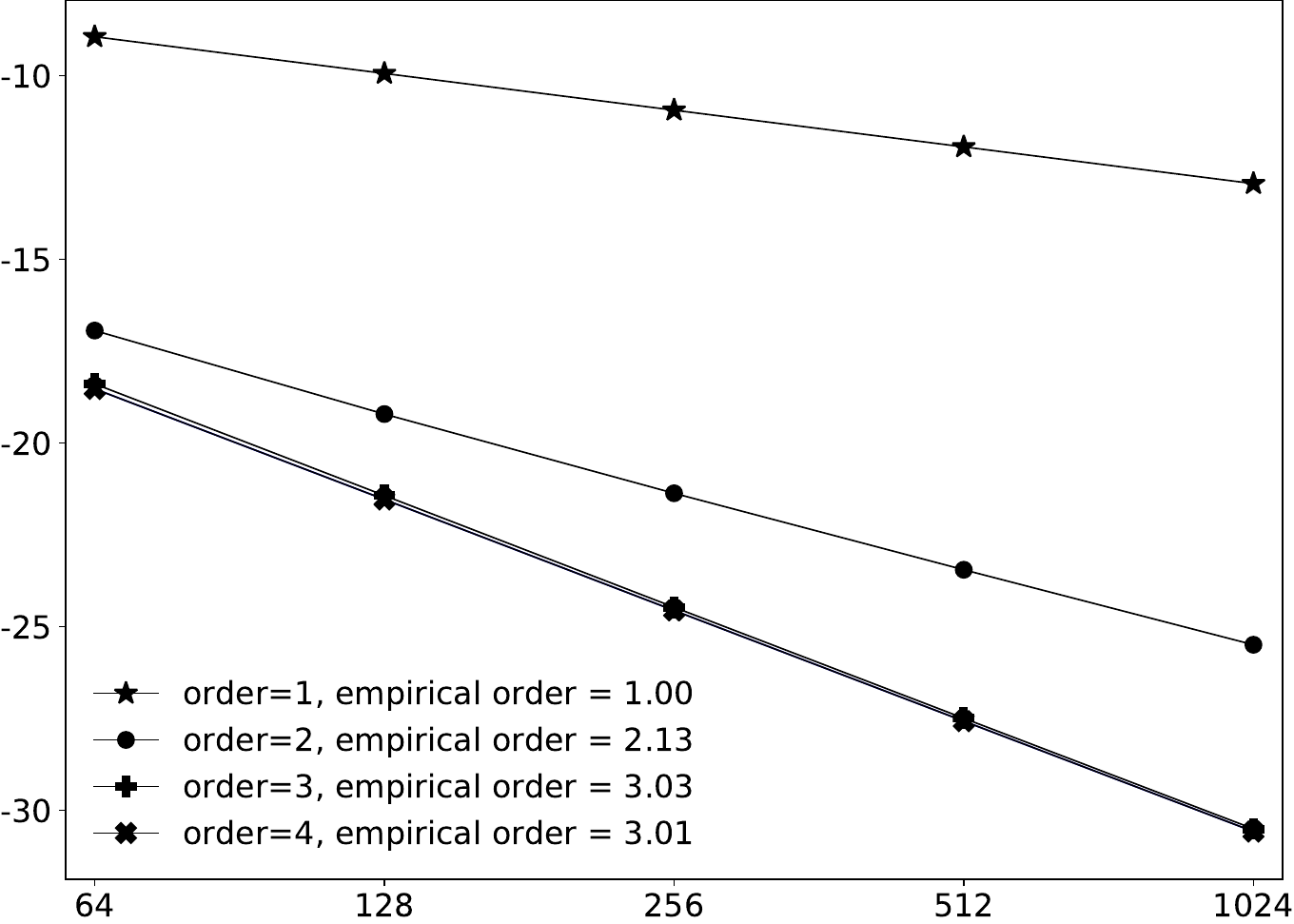}}
\caption{Same experiment as the one described in Figure \ref{fig-09}
with the initialization of the microscopic moments at first order following  (\ref{initialisation-ordre-1}).}
\label{fig-10} \end{figure} 

\renewcommand{\thetable}{2}
\begin{table}[H]
\smallskip
\centerline {\begin{tabular}{|c|c|c|c|c|}    \hline 
mesh points $ \backslash  $ equation order &  1  &  2 &  3 &  4 \\   \hline
initialization order  & 1 & 1 & 1 & 1  \\  \hline
  64 & 2.039 $10^{-3}$  & 7.967 $10^{-6}$  & 2.911 $10^{-6}$  & 2.652 $10^{-6}$   \\  \hline
 128 & 1.020 $10^{-3}$  & 1.648 $10^{-6}$  & 3.544 $10^{-7}$  & 3.290 $10^{-7}$   \\  \hline
 256 & 5.101 $10^{-4}$  & 3.697 $10^{-7}$  & 4.305 $10^{-8}$  & 4.049 $10^{-8}$   \\  \hline
 512 & 2.551 $10^{-4}$  & 8.730 $10^{-8}$  & 5.296 $10^{-9}$  & 5.018 $10^{-9}$   \\  \hline
1024 & 1.275 $10^{-4}$  & 2.120 $10^{-8}$  & 6.569 $10^{-10}$  & 6.247 $10^{-10}$   \\  \hline
convergence order & 1.00 & 2.13 & 3.03 & 3.01  \\  \hline
\end{tabular}}
\caption{Same numerical experiment as the one described in Table~\ref{tab-04}, except that
the initialization has  been changed to the first-order approximation (\ref{initialisation-ordre-1}). 
The precision is improved for second order and  we obtain the third order correctly, 
but the fourth-order approximation is only converging up to third order.}
\label{tab-05} \end{table}

\renewcommand{\thefigure}{7}
\begin{figure}    [h]  \centering
\centerline  {\includegraphics[width=.58\textwidth]   {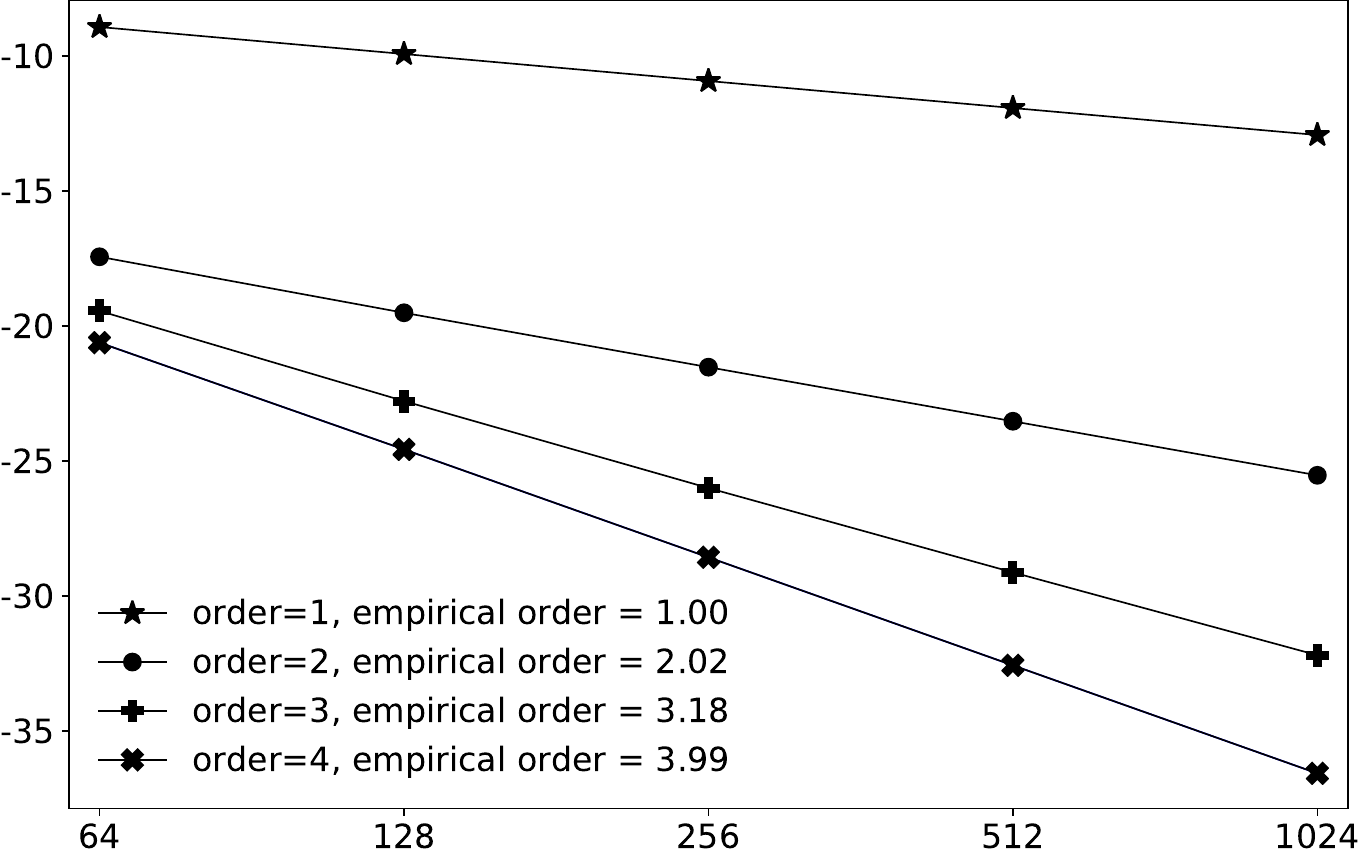}}
\vskip -.2 cm
\caption{Same experiment as the one described in Figure \ref{fig-09},  
with the initialization of the microscopic moments at second order following  (\ref{initialisation-ordre-2}).}
\label{fig-11} \end{figure} 

\renewcommand{\thetable}{3}
\begin{table}[H]
\centerline {\begin{tabular}{|c|c|c|c|c|}    \hline 
mesh points $ \backslash  $ equation order &  1  &  2 &  3 &  4 \\   \hline
initialization order  & 2 & 2 & 2 & 2  \\  \hline
  64 & 2.039 $10^{-3}$  & 5.607 $10^{-6}$  & 1.397 $10^{-6}$  & 6.191 $10^{-7}$   \\  \hline
 128 & 1.020 $10^{-3}$  & 1.332 $10^{-6}$  & 1.382 $10^{-7}$  & 3.997 $10^{-8}$   \\  \hline
 256 & 5.101 $10^{-4}$  & 3.299 $10^{-7}$  & 1.485 $10^{-8}$  & 2.506 $10^{-9}$   \\  \hline
 512 & 2.551 $10^{-4}$  & 8.233 $10^{-8}$  & 1.703 $10^{-9}$  & 1.567 $10^{-10}$   \\  \hline
1024 & 1.275 $10^{-4}$  & 2.057 $10^{-8}$  & 2.034 $10^{-10}$  & 9.798 $10^{-12}$   \\  \hline
convergence order & 1.00 & 2.02 & 3.18 & 3.99  \\  \hline
\end{tabular}}
\caption{Same numerical experiment as the one described in Table~\ref{tab-04}, except that
the initialization has  been changed to the second-order approximation (\ref{initialisation-ordre-2}). 
The precision order is now consistent with the approximation order.}
\label{tab-06} \end{table}

\renewcommand{\thetable}{4}
\begin{table}[H]
\centerline {\begin{tabular}{|c|c|c|c|c|}    \hline 
mesh points $ \backslash  $ equation order &  1  &  2 &  3 &  4 \\   \hline
initialization order  & 0 & 0 & 1 & 2  \\  \hline
  64 & 2.798 $10^{-3}$  & 7.606 $10^{-4}$  & 2.911 $10^{-6}$  & 6.191 $10^{-7}$   \\  \hline
 128 & 1.218 $10^{-3}$  & 1.983 $10^{-4}$  & 3.544 $10^{-7}$  & 3.997 $10^{-8}$   \\  \hline
 256 & 5.598 $10^{-4}$  & 4.979 $10^{-5}$  & 4.305 $10^{-8}$  & 2.506 $10^{-9}$   \\  \hline
 512 & 2.675 $10^{-4}$  & 1.245 $10^{-5}$  & 5.296 $10^{-9}$  & 1.567 $10^{-10}$   \\  \hline
1024 & 1.307 $10^{-4}$  & 3.113 $10^{-6}$  & 6.569 $10^{-10}$  & 9.798 $10^{-12}$   \\  \hline
convergence order & 1.10 & 1.99 & 3.03 & 3.99  \\  \hline
\end{tabular}}
\caption{Optimal initialization orders for the numerical experiment described in Table~\ref{tab-04}.
The precision order is now consistent with the approximation order without any extra calculation
for the initialization at the lowest  orders.}
\label{tab-07} \end{table}

\smallskip \noindent
In a second series of experiments with a constant velocity field, we use cubic parameters: 
\moneqstar
 U = 0.05 \,,\,\,   \alpha = -1 \,,\,\,  \sigma = 0.01 \,,
\monendstar
as previously. The second relaxation coefficient $ \, s' \, $
such that the relation (\ref{sigma-prime-cubique}) is satisfied:
\moneq  \label{s-prime-cubic}
\sigma' = 0.072425 \,,\,\,  s' =  1.7469537493994847 
\monend 
with $ \, \sigma' \equiv {1\over{s'}} - {1\over2} $.  
The initial condition is still a sine  wave
and we need only one term in the Fourier series.

\smallskip \noindent
In Table~\ref{tab-08}, the initialization for the second-order approximation is only of order zero and
the cubic convergence property is not obtained. On the other hand, when  the initialization for the  second-order
partial differential equation is first order accurate (see  Table~\ref{tab-09}), the second-order and third-order
approximations are identical.

\renewcommand{\thetable}{5}
\begin{table}[H]
\smallskip
\centerline {\begin{tabular}{|c|c|c|c|c|}    \hline 
mesh points $ \backslash  $ equation order &  1  &  2 &  3 &  4 \\   \hline
initialization order  & 0 & 0 & 1 & 2  \\  \hline
  64 & 2.826 $10^{-3}$  & 7.882 $10^{-4}$  & 2.181 $10^{-6}$  & 6.455 $10^{-7}$   \\  \hline
 128 & 1.219 $10^{-3}$  & 1.988 $10^{-4}$  & 2.332 $10^{-7}$  & 4.063 $10^{-8}$   \\  \hline
 256 & 5.598 $10^{-4}$  & 4.980 $10^{-5}$  & 2.665 $10^{-8}$  & 2.544 $10^{-9}$   \\  \hline
 512 & 2.675 $10^{-4}$  & 1.245 $10^{-5}$  & 3.174 $10^{-9}$  & 1.590 $10^{-10}$   \\  \hline
1024 & 1.307 $10^{-4}$  & 3.113 $10^{-6}$  & 3.869 $10^{-10}$  & 9.951 $10^{-12}$   \\  \hline
convergence order & 1.11 & 2.00 & 3.11 & 4.00  \\  \hline
\end{tabular}}
\caption{Errors measured with the maximum norm
  between the D1Q3 lattice Boltzmann  scheme with cubic parameter 
(\ref{s-prime-cubic}) for   
$ \, \alpha = -1 $,  $ \, \sigma \equiv {{1}\over{s}} - {1\over2} = 0.01 $. 
Even if the relaxation parameters have been fitted in order to obtain third-order accuracy
with the second-order equivalent partial differential equation, 
the error remains second-order accurate in this case.}
\label{tab-08} \end{table}

\renewcommand{\thetable}{6}
\begin{table}[H]
\smallskip
\centerline {\begin{tabular}{|c|c|c|c|c|}    \hline 
mesh points $ \backslash  $ equation order &  1  &  2 &  3 &  4 \\   \hline
initialization order  & 0 & 1 & 1 & 2  \\  \hline
  64 & 2.826 $10^{-3}$  & 2.181 $10^{-6}$  & 2.181 $10^{-6}$  & 6.455 $10^{-7}$   \\  \hline
 128 & 1.219 $10^{-3}$  & 2.332 $10^{-7}$  & 2.332 $10^{-7}$  & 4.063 $10^{-8}$   \\  \hline
 256 & 5.598 $10^{-4}$  & 2.665 $10^{-8}$  & 2.665 $10^{-8}$  & 2.544 $10^{-9}$   \\  \hline
 512 & 2.675 $10^{-4}$  & 3.174 $10^{-9}$  & 3.174 $10^{-9}$  & 1.590 $10^{-10}$   \\  \hline
1024 & 1.307 $10^{-4}$  & 3.869 $10^{-10}$  & 3.869 $10^{-10}$  & 9.951 $10^{-12}$   \\  \hline
convergence order & 1.11 & 3.11 & 3.11 & 4.00  \\  \hline
\end{tabular}}
\caption{Same numerical experiment as the one described in Table~\ref{tab-08},
except that the initialization scheme is first-order accurate when comparing with
the second-order equivalent partial differential equation. The third-order terms
of the partial differential equation are identically null in this case due to the choice of a set of cubic parameters,
and the order of accuracy jumps to third order.}
\label{tab-09} \end{table}

\smallskip \noindent
When the velocity is no longer constant but given by the relation (\ref{u-cosinus}),
the modes are coupled as detailed in Section 5. 
We have used 30 active modes in the Fourier series. 
In Tables~\ref{tab-10} to \ref{tab-12}, we experiment with the three types of initialization, 
(\ref{initialisation-ordre-0}), (\ref{initialisation-ordre-1}) and (\ref{initialisation-ordre-2}).
The results are qualitatively identical to the previous experiments with a uniform vector field.
When the initialization is done with the equilibrium (\ref{initialisation-ordre-0}), the lattice Boltzmann
scheme can  be compared with equivalent partial  differential equations only at second order, 
as detailed in Table~\ref{tab-10}.
For the first order  (\ref{initialisation-ordre-1}), third-order accuracy can be obtained.
Nevertheles, the fourth-order differential model is only third-order accurate (Table~\ref{tab-11}). 
With a  second-order initialization~(\ref{initialisation-ordre-2}), the asymptotic partial
differential equation of a given degree is an approximation of the lattice Boltzmann scheme 
with the same degree, as presented in Table~\ref{tab-12}.

\renewcommand{\thetable}{7}
\begin{table}[H]
\smallskip
\centerline {\begin{tabular}{|c|c|c|c|c|}    \hline 
mesh points $ \backslash  $ equation order &  1  &  2 &  3 &  4 \\   \hline
initialization order  & 0 & 0 & 0 & 0  \\  \hline
  64 & 5.625 $10^{-3}$  & 1.172 $10^{-3}$  & 1.120 $10^{-3}$  & 1.088 $10^{-3}$   \\  \hline
 128 & 2.534 $10^{-3}$  & 2.952 $10^{-4}$  & 2.819 $10^{-4}$  & 2.778 $10^{-4}$   \\  \hline
 256 & 1.195 $10^{-3}$  & 7.327 $10^{-5}$  & 6.992 $10^{-5}$  & 6.940 $10^{-5}$   \\  \hline
 512 & 5.793 $10^{-4}$  & 1.823 $10^{-5}$  & 1.739 $10^{-5}$  & 1.733 $10^{-5}$   \\  \hline
1024 & 2.851 $10^{-4}$  & 4.547 $10^{-6}$  & 4.337 $10^{-6}$  & 4.329 $10^{-6}$   \\  \hline
convergence order & 1.07 & 2.00 & 2.00 & 2.00  \\  \hline
\end{tabular}}
\caption{Same numerical experiment as the one described in Table~\ref{tab-04};
the uniform vector field is replaced by a cosine velocity (\ref{u-cosinus}) with $ \, U_0 = 0.05$. 
As in the previous experiment, the error remains second-order accurate even if we use the third-order or the fourth-order
equivalent  equation for the approximation of the lattice Boltzmann scheme.}
\label{tab-10} \end{table}

\renewcommand{\thetable}{8}
\begin{table}[H]
\smallskip
\centerline {\begin{tabular}{|c|c|c|c|c|}    \hline 
mesh points $ \backslash  $ equation order &  1  &  2 &  3 &  4 \\   \hline
initialization order  & 1 & 1 & 1 & 1  \\  \hline
  64 & 4.561 $10^{-3}$  & 1.087 $10^{-4}$  & 5.628 $10^{-5}$  & 2.446 $10^{-5}$   \\  \hline
 128 & 2.259 $10^{-3}$  & 1.982 $10^{-5}$  & 6.495 $10^{-6}$  & 2.427 $10^{-6}$   \\  \hline
 256 & 1.126 $10^{-3}$  & 4.124 $10^{-6}$  & 7.747 $10^{-7}$  & 2.622 $10^{-7}$   \\  \hline
 512 & 5.620 $10^{-4}$  & 9.336 $10^{-7}$  & 9.441 $10^{-8}$  & 3.018 $10^{-8}$   \\  \hline
1024 & 2.808 $10^{-4}$  & 2.216 $10^{-7}$  & 1.165 $10^{-8}$  & 3.610 $10^{-9}$   \\  \hline
convergence order & 1.01 & 2.23 & 3.06 & 3.18  \\  \hline
\end{tabular}}
\caption{Same numerical experiment as the one described in Table~\ref{tab-10}. 
The initialization is now given by the first-order approximation (\ref{initialisation-ordre-1}). 
The precision is improved for the second-order partial differential equation and  
we obtain the third-order correctly.
But the fourth-order approximation is converging only up to third order.}
\label{tab-11} \end{table}

\renewcommand{\thetable}{9}
\begin{table}[H]
\smallskip
\centerline {\begin{tabular}{|c|c|c|c|c|}    \hline 
mesh points $ \backslash  $ equation order &  1  &  2 &  3 &  4 \\   \hline
initialization order  & 2 & 2 & 2 & 2  \\  \hline
  64 & 4.548 $10^{-3}$  & 9.505 $10^{-5}$  & 4.264 $10^{-5}$  & 1.082 $10^{-5}$   \\  \hline
 128 & 2.257 $10^{-3}$  & 1.807 $10^{-5}$  & 4.742 $10^{-6}$  & 6.741 $10^{-7}$   \\  \hline
 256 & 1.125 $10^{-3}$  & 3.904 $10^{-6}$  & 5.540 $10^{-7}$  & 4.162 $10^{-8}$   \\  \hline
 512 & 5.620 $10^{-4}$  & 9.060 $10^{-7}$  & 6.678 $10^{-8}$  & 2.544 $10^{-9}$   \\  \hline
1024 & 2.808 $10^{-4}$  & 2.182 $10^{-7}$  & 8.191 $10^{-9}$  & 1.522 $10^{-10}$   \\  \hline
convergence order & 1.00 & 2.19 & 3.08 & 4.03  \\  \hline
\end{tabular}}
\caption{Same numerical experiment as the one described in Table~\ref{tab-10}. 
  The initialization is now given by the second order approximation (\ref{initialisation-ordre-2}).
The precision order is now consistent  with the approximation order.}
\label{tab-12} \end{table}

\smallskip \noindent
  Optimal initialization orders for the numerical experiment with sinusoidal velocity
  can be made precise  as follows:
\moneqstar  \left\{ \begin{array} {rll}
{\textrm {partial differential equation order}} &= 1 {\textrm { or }} 2:  & {\textrm {initialization at order 0}} \\
{\textrm {partial differential equation order}} &= 3: & {\textrm {initialization at order 1}} \\
{\textrm {partial differential equation order}} &= 4: & {\textrm {initialization at order 2}}.
\end{array} \right. \monendstar
The precision order is now consistent  with the approximation order without any extra  calculation
for the initialization at the lowest  orders.

\smallskip \noindent
If  the initial condition is no longer  a  sinusoidal  wave but a constant state, the results presented in
Tables~\ref{tab-10} to \ref{tab-12} are essentially unchanged.
We present in Table~\ref{tab-13} the analogue  of Table~\ref{tab-08} for this case. 

\renewcommand{\thetable}{10}
\begin{table}[H]
\smallskip
\centerline {\begin{tabular}{|c|c|c|c|c|}    \hline
mesh points $ \backslash  $ equation order &  1  &  2 &  3 &  4 \\   \hline
initialization order  & 0 & 0 & 1 & 2  \\  \hline
  64 & 6.050 $10^{-4}$  & 3.597 $10^{-5}$  & 1.224 $10^{-5}$  & 1.306 $10^{-6}$   \\  \hline
 128 & 2.932 $10^{-4}$  & 7.528 $10^{-6}$  & 1.475 $10^{-6}$  & 8.102 $10^{-8}$   \\  \hline
 256 & 1.447 $10^{-4}$  & 1.699 $10^{-6}$  & 1.800 $10^{-7}$  & 5.034 $10^{-9}$   \\  \hline
 512 & 7.194 $10^{-5}$  & 4.024 $10^{-7}$  & 2.221 $10^{-8}$  & 3.130 $10^{-10}$   \\  \hline
1024 & 3.587 $10^{-5}$  & 9.784 $10^{-8}$  & 2.758 $10^{-9}$  & 1.943 $10^{-11}$   \\  \hline
convergence order & 1.02 & 2.13 & 3.03 & 4.01  \\  \hline
\end{tabular}}
\caption{Optimal initialization orders for the numerical experiment with sinusoidal velocity
  described in Table~\ref{tab-10}. The initial condition is changed from a sinusoidal function
  to a constant state. 
  Each asymptotic partial differential
  equation presents a precision order consistent with its approximation order.}
\label{tab-13} \end{table}

\smallskip \noindent 
We tried also to apply a cubic choice of coefficients for the non homogeneous case.
We have not observed any spectacular improved precision.
There is no inconsistency because the cubic parameters have been explicated with the hypothesis
of a constant velocity field.

\bigskip \bigskip    \noindent {\bf \large    9) \quad  Long-time asymptotic study}

\smallskip \noindent 
It is always difﬁcult to reconcile an expansion of a  regular solution when the solution is steady.
For instance, in steady forced Poiseuille ﬂow the body force must balance the viscous stress.
However, the body force typically appears at leading order, 
while the viscous stress arises at 
first order relative to the spatial step. 
Then whithin 
the acoustic scaling framework, 
we obtain different stationary solutions for different meshes as the number of
mesh points increases.  
We refer, {\it e.g.},
to the contribution 
\cite{DLT20} for  Poiseuille flow with anti-bounce-back  
boundary conditions.

\smallskip \noindent
In this section, we first present the numerical analysis of the spectral properties
of the D1Q3 scheme. Then we derive a simple numerical method to achieve a stationary solution at second order.
Finally the spectral approach developed in the previous sections is adapted to this stationary case.
We report the results of the  
numerical simulations.
We conclude with a preliminary  
analysis of our analytical and numerical results.

\smallskip \noindent
We put 
in evidence some intrinsic properties of the D1Q3 lattice Boltzmann scheme
with the first unstationary mode. One step of the algorithm on a grid with $ \, N \, $ mesh points can be written
\moneq    \label{iterations-schema}
f(t+\Delta t) = A_{\textrm D1Q3} \,\, f(t)
\monend 
with $ \, A_{\textrm D1Q3} \, $  the global iteration matrix of order $ \, 3 N \times 3 N \, $
of this linear scheme. The matrix  $ \, A_{\textrm D1Q3} \, $  contains all  information
relative to   collision and  advection for all the vertices. 
With an Arnoldi algorithm  (see {\it e.g.} \cite{LL00}), we extract the first eigenmode of the matrix $ \,  A_{\textrm D1Q3} $.
This eigenvalue~$ \, \gamma \, $ 
is numerically real in our case and we introduce a scaled  parameter~$\, \Gamma \, $
defined as follows. 
From the operator $ \, \alpha_2 \, $ in (\ref{alpha2-beta2}), we first introduce the discrete equivalent diffusivity
$ \,\, \kappa = \lambda \, \Delta x \, \sigma \, {{\alpha+2}\over3} $. Then for a simulation with a wave number
$ \, k $, we set
\moneq    \label{Gamma_mode}
\Gamma = -{{\gamma}\over{\kappa \, k^2}} 
\monend
and the minus sign is introduced for positive numbers.
With this definition of $\,  \Gamma $,  the diffusion equation leads to a numerical value of $  \Gamma $ equal some integer. 
This is the scaled  first eigenvalue of the iteration matrix  $ \, A_{\textrm D1Q3} $.
Then from the corresponding eigenvector $\,  f_{\gamma} $, we extract the conserved moment
\moneq    \label{densite-pierre}
\rho_\gamma  = \sum_{j=1}^{j=3}   f_{\gamma,\, j}  .
\monend
It is a function defined at all mesh points. 
We have represented in   Figures \ref{modes-a} to  \ref{modes-d} the corresponding modes for
$ \, U = 0 $, $ \, 0.0005$, $ \, 0.005 \, $ and $ \, 0.05 $.

\smallskip \noindent 
The very interesting observation concerns the evolution of the eigenvalue $ \, \Gamma \, $
as  function of velocity and  number of mesh points
($N = 64$, $\, 128$, $ \, 356 $, $\, 512 $) presented in Figure \ref{modes-Gamma}. 
A spectacular growth occurs for the largest velocity. In practice, the lattice Boltzmann scheme is much more viscous
than proposed by the natural scaling $ \,\, \kappa \, k^2 $.
%
We have also observed that for large  values of the velocity $ \, U $,
we find values for $ \, \Gamma \, $ roughly proportional to $ \, U \, N $.

\smallskip \noindent 
We have done numerical experiments with three advective velocities given by the relation~(\ref{u-cosinus}) with 
$ \, U = 0.0005 $, $ \, U = 0.005 \, $ and  $ \, U = 0.05 $.
For each of these parameters, we have used  four discretizations with  64, 128, 256 and 
512 mesh points.
We have made various choices for the approximation of the D1Q3 stationary field.

\renewcommand{\thefigure}{8}
\begin{figure}    [H]  \centering
\centerline  {\includegraphics[width=.58\textwidth]   {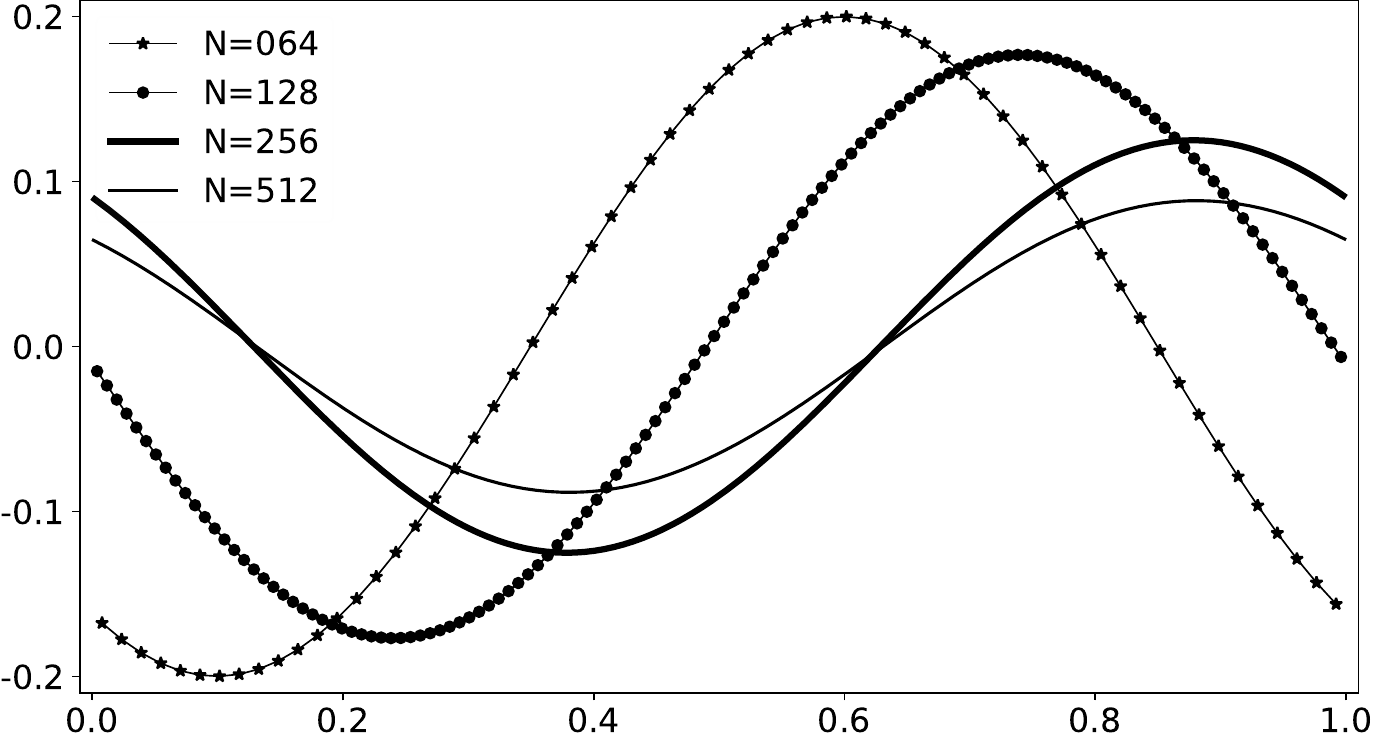}}
  
\caption{First eigenmode of the stationary D1Q3 discrete dynamics, $ \, U = 0  $.
We have
$ \, \Gamma_{64}=1.00053560 $, 
$ \, \Gamma_{128}=1.00013382 $, 
$ \, \Gamma_{256}=1.00003341 \, $ and 
    $ \, \Gamma_{512}=1.00000831 $.   
The difference between $\Gamma$ and its integer part is due to rounding errors.}
\label{modes-a} \end{figure} 

\renewcommand{\thefigure}{9}
\begin{figure}    [H]  \centering
  \centerline  {\includegraphics[width=.58\textwidth]   {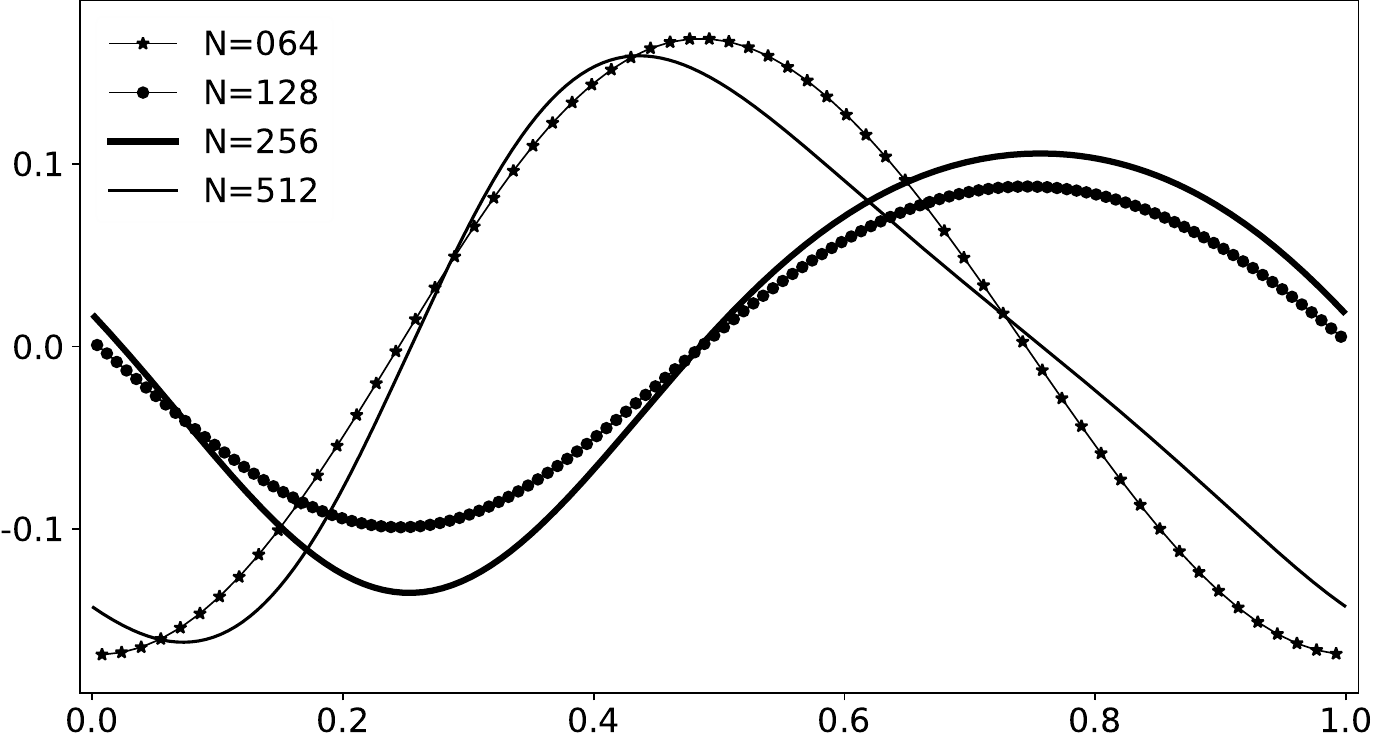}}
\caption{First eigenmode of the stationary D1Q3 discrete dynamics, $ \, U = 0.0005  $.
We have
$ \, \Gamma_{64}=1.00193081 $, 
$ \, \Gamma_{128}=1.00572981 $, 
$ \, \Gamma_{256}=1.02241778 \, $ and 
$ \, \Gamma_{512}=1.08930892 $.} 
\label{modes-b} \end{figure} 

\renewcommand{\thefigure}{10}
\begin{figure}    [H]  \centering
\centerline  {\includegraphics[width=.58\textwidth]   {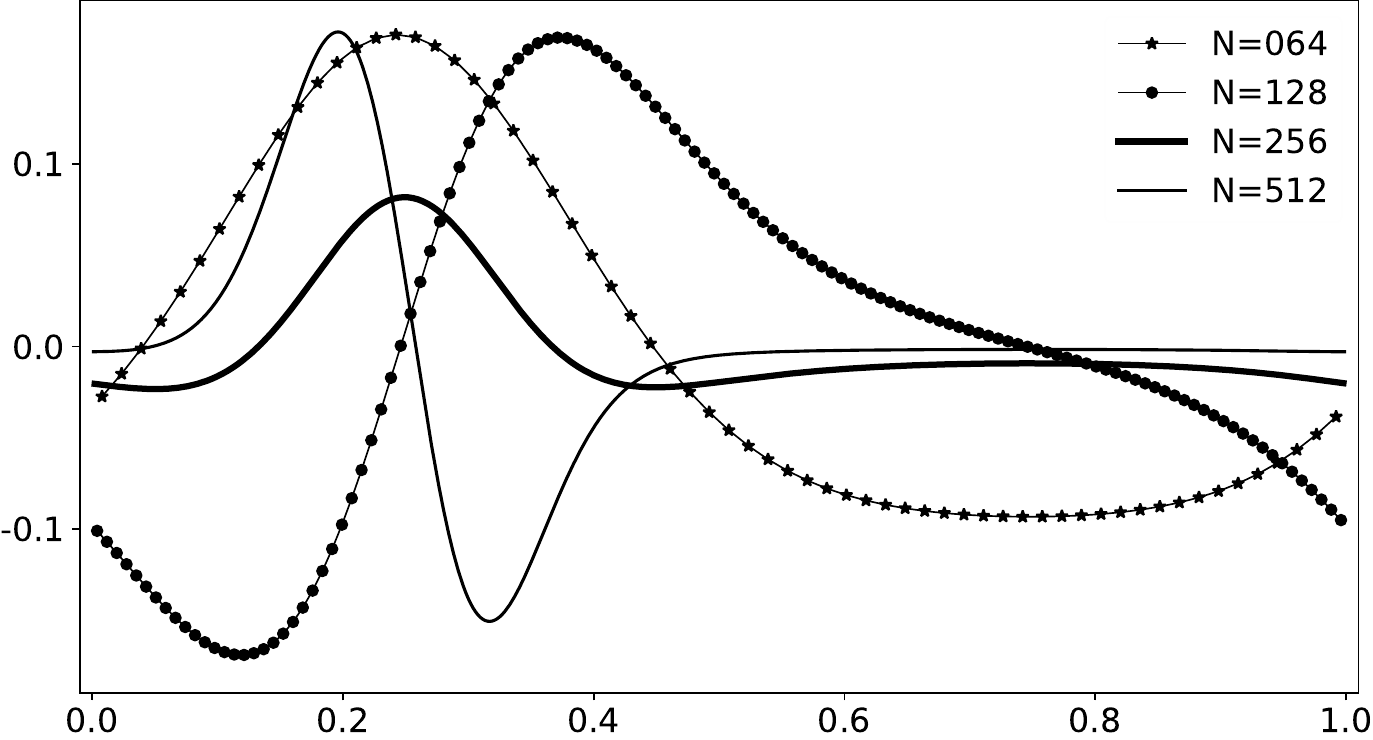}}
\caption{First eigenmode of the stationary D1Q3 discrete dynamics, $ \, U = 0.005  $.
We have
$ \, \Gamma_{64}=1.13928618 $, 
$ \, \Gamma_{128}=1.54649116 $, 
$ \, \Gamma_{256}=3.00825658  \, $ and 
$ \, \Gamma_{512}=6.75648764 $.} 
\label{modes-c} \end{figure} 

\renewcommand{\thefigure}{11}
\begin{figure}    [H]  \centering
\centerline  {\includegraphics[width=.58\textwidth]   {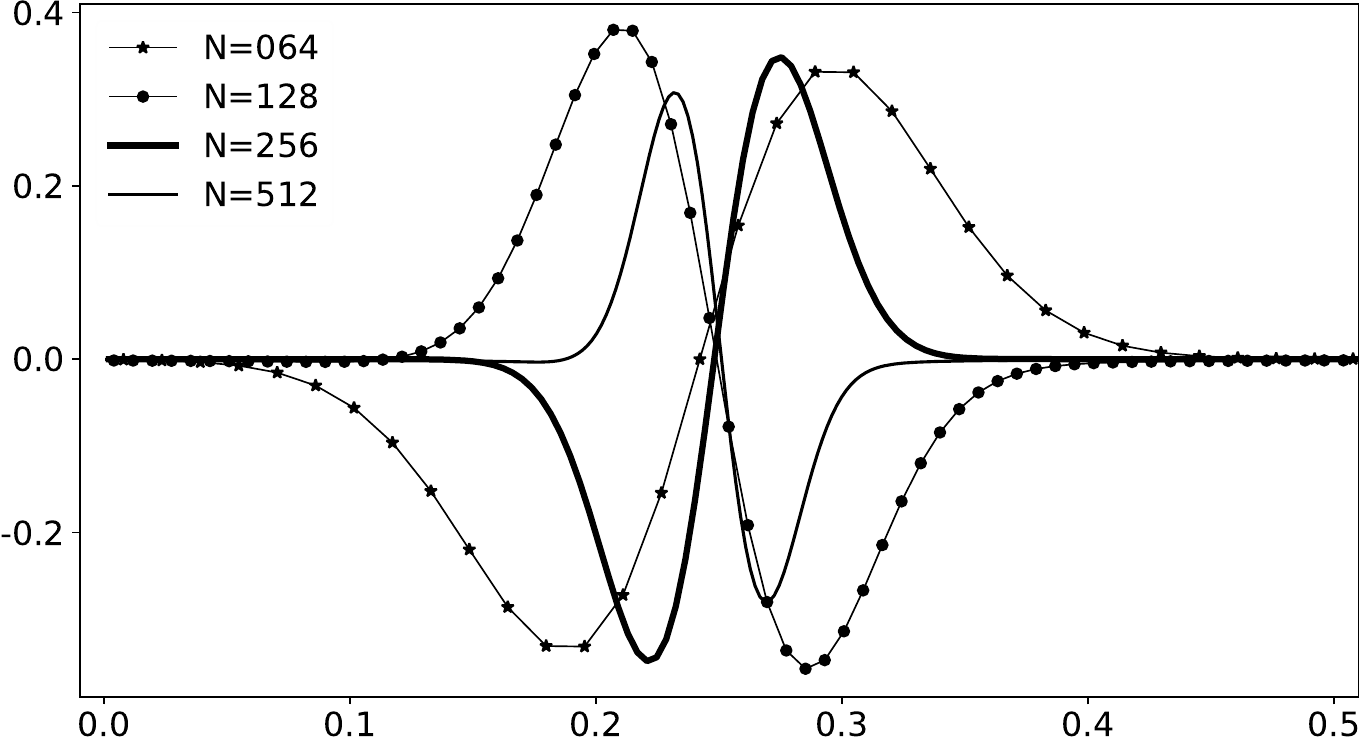}}
\caption{First eigenmode of the stationary D1Q3 discrete dynamics, $ \, U = 0.05  $.
We have
$ \, \Gamma_{64}=8.62260312 $, 
$ \, \Gamma_{128}=17.81972445 $, 
$ \, \Gamma_{256}=36.16622885   \, $ and 
$ \, \Gamma_{512}=72.84095421 $.} 

\label{modes-d} \end{figure} 

\renewcommand{\thefigure}{12}
\begin{figure}    [H]  \centering
\centerline  {\includegraphics[width=.58 \textwidth]   {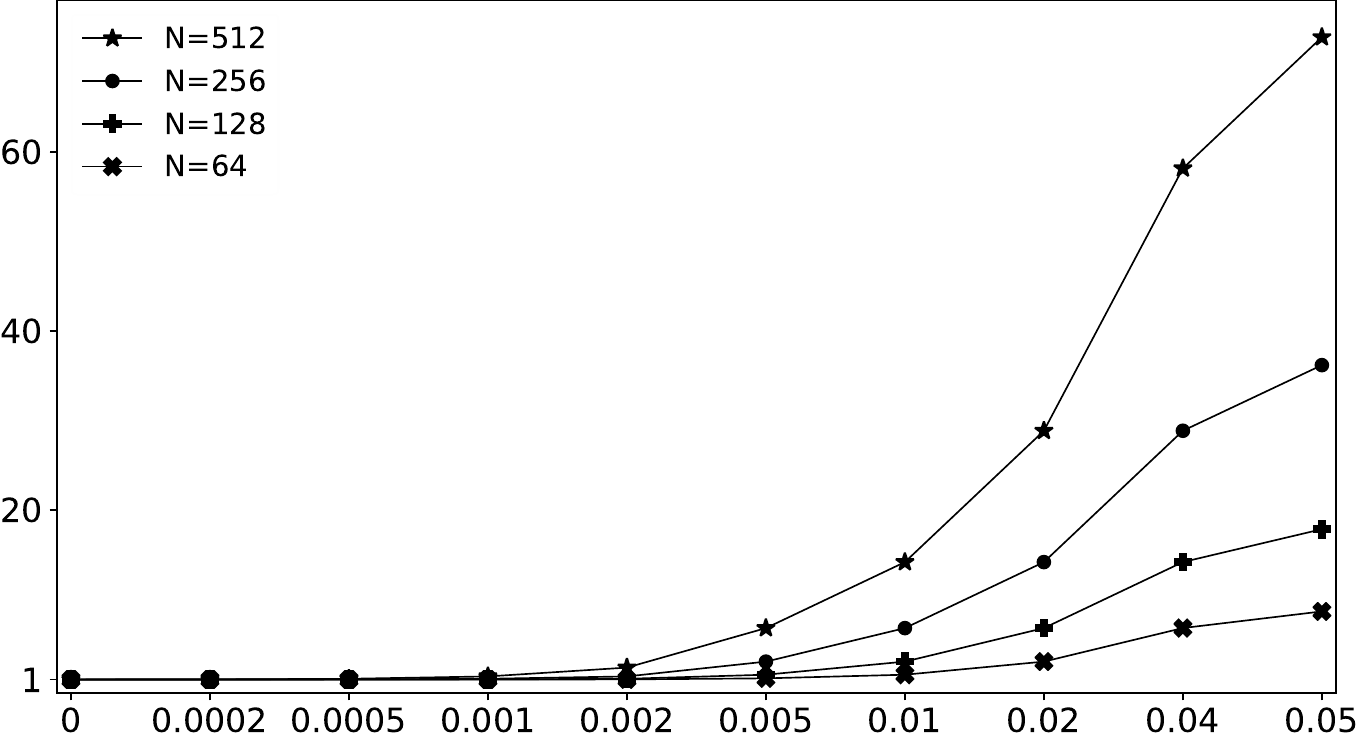}}
\caption{First eigenvalue $ \, \Gamma \, $ defined in (\ref{Gamma_mode})  for the stationary D1Q3 discrete dynamics
as a function of velocity and number of mesh points. Observe that the $x$-scale is neither linear nor logarithmic
to clearly highlight the numerical values.}
\label{modes-Gamma} \end{figure} 

\smallskip \noindent
When we study the equivalent partial differential equation of the lattice Boltzmann scheme, 
we remark that the equation 
\moneqstar 
\partial_t \rho + \partial_x \big( \lambda\, U\, \cos(k \, x) \, \rho \big) - \, \mu \, \partial^2_x \rho = 0
\monendstar
is a simple approximation at second-order accuracy. 
Moreover, we 
observe that the  analytical expression of the stationary solution  for the equation
with the integral condition
\moneq    \label{integrale-unite}
\int_0^L \rho(x) \, \dd x = 1  . 
\monend 
can be made explicit  as 
\moneq    \label{stationnaire-analytique}
\rho(x) = K \, \exp \Big( {{\lambda\,U}\over{k\, \mu}}\, \sin (k \, x) \Big) .
\monend
The normalization constant $ \, K \, $ in relation (\ref{stationnaire-analytique})
is chosen such that the condition (\ref{integrale-unite}) is satisfied.

\smallskip \noindent
We then 
compare the numerical solution obtained with the D1Q3 scheme with the numerical solution
of Fourier series truncated with  30  active  modes.
%
We introduce an operator $ \, A \, $ obtained at various orders from the relation (\ref{operateur-A}) typically and
\moneqstar
A_\infty =  {{1}\over{\partial_x}} \, A . 
\monendstar
Then $ \, A_\infty \, \rho = {\textrm{constant}} \, $ and this constant is zero by periodicity
of all the  functions of the problem.
Then $ \, A_\infty \, \rho = 0 \, $ with an operator  $ \, A_\infty \, $ given at various orders 
by  
\moneqstar
A_\infty^1 = \lambda \, m_u   -\mu \, \partial_x
\monendstar
at order 1,
\moneqstar
A_\infty^2 = A_\infty^1  + \mu \, m_u \, \partial_x \, m_u 
\monendstar
at order 2 and 
\moneqstar
A_\infty^3 = A_\infty^2  + \xi_u \, m_u \, \partial_x \, m_u \, \partial_x \, m_u
+ \xi_{xu}  \,\partial_x^2 \,  m_u   + \xi_{ux}  \, m_u \, \partial_x^2 
\monendstar
at order 3 and finally, 
\moneqstar \left\{ \begin{array} {l}
A_\infty^4 = A_\infty^3   
+  \zeta_{u4} \, m_u \, \partial_x \, m_u \, \partial_x \, m_u \, \partial_x \, m_u 
+ \zeta_{xxuu}  \, \partial_x^2 \,  m_u \, \partial_x \, m_u  +  \zeta_{uxxu}  \, m_u \, \partial_x^3 \, \, m_u \\ 
\qquad \qquad + \zeta_{uuxx}  \, m_u \, \partial_x \, m_u  \, \partial_x^2 +  \zeta_{x4} \, \partial_x^3 
\end{array} \right. \monendstar
at order 4.

\smallskip \noindent
For $ \, U = 0.0005 $, 
The parameters  of the D1Q3 scheme are identical to those chosen for unstationary  simulations 
[$\alpha = -1 $, $ \, s=1.5 $, $ \, s' = 1.2 $] and we again use 30 Fourier modes. 
the quantitative results  are presented in Table~\ref{tab-01}.
  With $\, s'=1.5 \, $ the convergence results are very similar to the ones of Table~\ref{tab-01}.

\renewcommand{\thetable}{11}
\begin{table}[H]
\smallskip
\centerline {\begin{tabular}{|c|c|c|c|c|}    \hline 
mesh points $ \backslash  $ equation order & 1  &  2 &  3 &  4 \\   \hline 
$ 64  $&$  8.182 \, 10^{-5} $&$    8.167 \, 10^{-5} $&$  5.455 \, 10^{-5} $&$  6.935  \, 10^{-8} $ \\   \hline 
$ 128 $&$  4.495 \, 10^{-5} $&$    4.483 \, 10^{-5} $&$  2.997 \, 10^{-5} $&$  1.113  \, 10^{-8} $ \\   \hline 
$ 256 $&$  2.616 \, 10^{-5} $&$    2.611 \, 10^{-5} $&$  1.744 \, 10^{-5} $&$  2.067  \, 10^{-9} $ \\   \hline 
$ 512 $&$  1.601 \, 10^{-5} $&$    1.610 \, 10^{-5} $&$  1.068 \, 10^{-5} $&$  4.836  \, 10^{-10} $ \\   \hline 
convergence order &$ 0.78 $&$  0.78  $&$0.78$&$ 2.39   $ \\   \hline 
\end{tabular}}

\caption{Differences between the lattice Boltzmann D1Q3 scheme and various equivalent equations
for a stationary experiment with $ \, U = 0.0005 $.}
\label{tab-01} \end{table}

\renewcommand{\thetable}{12}
\begin{table}[H]
\smallskip
\centerline {\begin{tabular}{|c|c|c|c|c|}    \hline 
mesh points $ \backslash  $ equation order & 1  &  2 &  3 &  4 \\   \hline 
$ 64  $&$  1.362 \, 10^{-3} $&$    1.378 \, 10^{-3} $&$  9.083 \, 10^{-4} $&$  2.886  \, 10^{-6} $ \\   \hline 
$ 128 $&$  8.538 \, 10^{-4} $&$    8.845 \, 10^{-4} $&$  5.692 \, 10^{-4} $&$  8.780  \, 10^{-7} $ \\   \hline 
$ 256 $&$  6.183 \, 10^{-4} $&$    6.437 \, 10^{-4} $&$  4.122 \, 10^{-4} $&$  3.066  \, 10^{-7} $ \\   \hline 
$ 512 $&$  4.578 \, 10^{-4} $&$    4.750 \, 10^{-4} $&$  3.052 \, 10^{-4} $&$  1.052  \, 10^{-7} $ \\   \hline 
convergence order &$ 0.52 $&$ 0.51 $&$ 0.52  $&$  1.58  $ \\   \hline 
\end{tabular}}
\caption{Differences between the lattice Boltzmann D1Q3 scheme and various equivalent equations
for a stationary experiment with $ \, U = 0.005 $.}
\label{tab-02} \end{table}

\renewcommand{\thetable}{13}
\begin{table}[H]
\smallskip
\centerline {\begin{tabular}{|c|c|c|c|c|}    \hline 
mesh points $ \backslash  $ equation order & 1  &  2 &  3 &  4 \\   \hline 
$ 64  $&$  3.883 \, 10^{-2} $&$    4.042 \, 10^{-2} $&$  2.590 \, 10^{-2} $&$  6.585  \, 10^{-4} $ \\   \hline 
$ 128 $&$  2.856 \, 10^{-2} $&$    2.967 \, 10^{-2} $&$  1.904 \, 10^{-2} $&$  2.439  \, 10^{-4} $ \\   \hline 
$ 256 $&$  2.057 \, 10^{-2} $&$    2.136 \, 10^{-2} $&$  1.372 \, 10^{-2} $&$  8.820  \, 10^{-5} $ \\   \hline 
$ 512 $&$  1.468 \, 10^{-2} $&$    1.523 \, 10^{-2} $&$  9.790 \, 10^{-3} $&$  3.153  \, 10^{-5} $ \\   \hline 
convergence order &$ 0.47  $&$  0.47  $&$  0.47  $&$  1.46 $ \\   \hline 
\end{tabular}}
\caption{Differences between the lattice Boltzmann D1Q3 scheme and various equivalent equations
for a stationary experiment with $ \, U = 0.05 $.}
\label{tab-03} \end{table}


\smallskip \noindent
We observe that increasing   the  order of accuracy increases 
the quality of the approximation
between the  lattice Boltzmann scheme and the computation with  Fourier series.
We observe that 
 it takes very many timesteps before the solution approaches its large-time limit to sufﬁcient accuracy.
This is consistent with the relaxation diffusion time $ \, \tau = {{\lambda}\over{\mu \, \Delta x \, k^2}} \, $
measured with our scaling.  
%
For example, with  512 mesh points, we have used more than\br
3,400,000 time steps to reach  the numerical result presented in Table~\ref{tab-03}.
This difficulty in reaching the stationary state
is  directly correlated with the high value
  of the eigenvalue~(72.8) in this case. 
We observe also that the stationary solution given by the 
asymptotic expansion is globally correct with errors between $ \, 10^{-5} \, $ and  $ \, 10^{-10} $.
The convergence order is, however, 
slower than expected with the order of the
partial differential equation. Nevertheless, 
the convergence order for the formal fourth-order approximation (fourth column of
Table~\ref{tab-01}) is~2.39.

\smallskip \noindent 
For $ \, U = 0.005 $,
we present the numerical results in Figures \ref{fig-04-a} and  \ref{fig-04-b} 
and the quantitative residuals in Table~\ref{tab-02}.
The parameters  of the D1Q3 scheme are identical to those chosen for the other simulations.
The convergence process when the mesh is refined is slow.
We have, for example, for the fourth-order partial differential equation
(fourth column in Table~\ref{tab-02}), that  a least-square fitting gives a convergence order of 1.58. 

\renewcommand{\thefigure}{13}
\begin{figure}    [H]  \centering
\centerline  {\includegraphics[width=.58 \textwidth]   {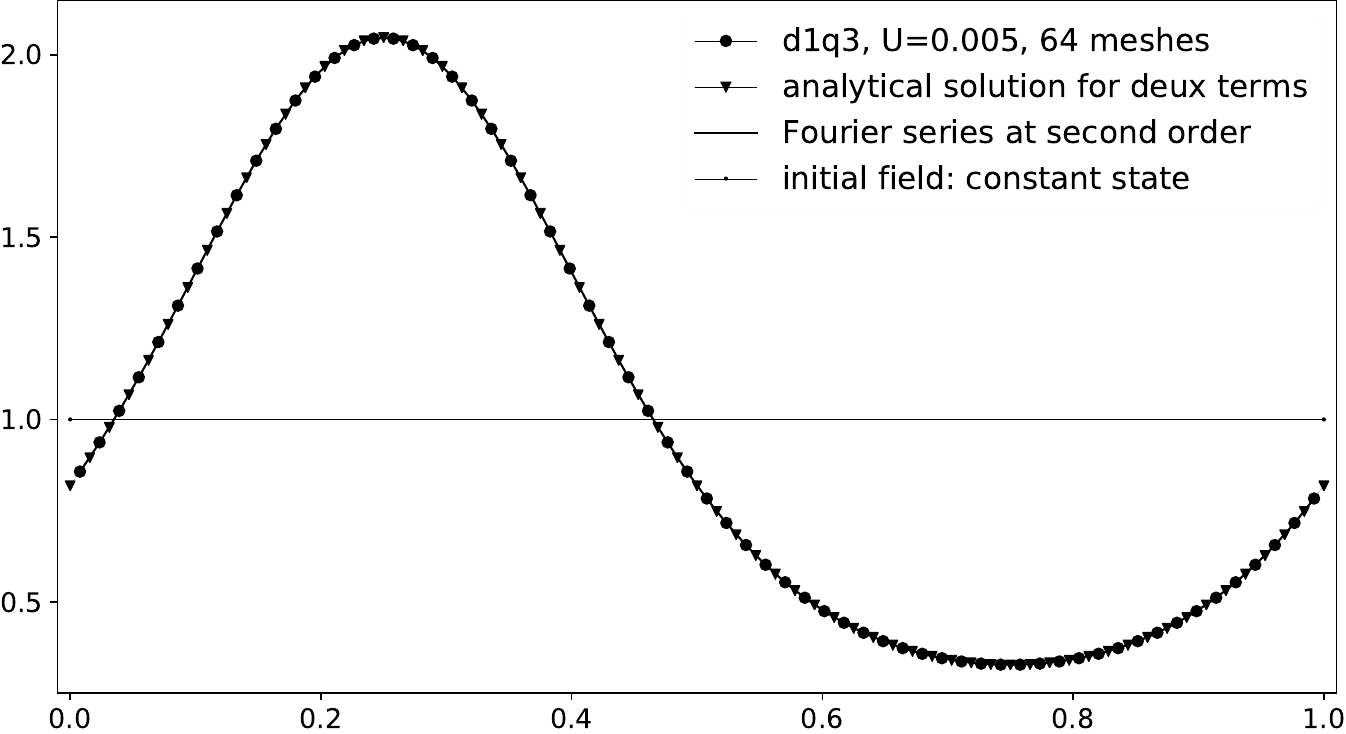}}
\caption{Stationary field for a sinusoidal advection field with  $ \, U = 0.005  $ and  64 mesh points.
 The analytic formula (\ref{stationnaire-analytique}) obtained
with the advective terms and only the uniform dissipation gives a very correct approximation
of the stationary asymptotic solution obtained with the D1Q3 lattice Boltzmann scheme.
It is just necessary to compute precisely the constant $ \, K \, $ in order to satisfy the  integral condition
(\ref{integrale-unite}).  
}  
\label{fig-04-a} \end{figure} 

\renewcommand{\thefigure}{14}
\begin{figure}    [H]  \centering
  \centerline  {\includegraphics[width=.58\textwidth]   {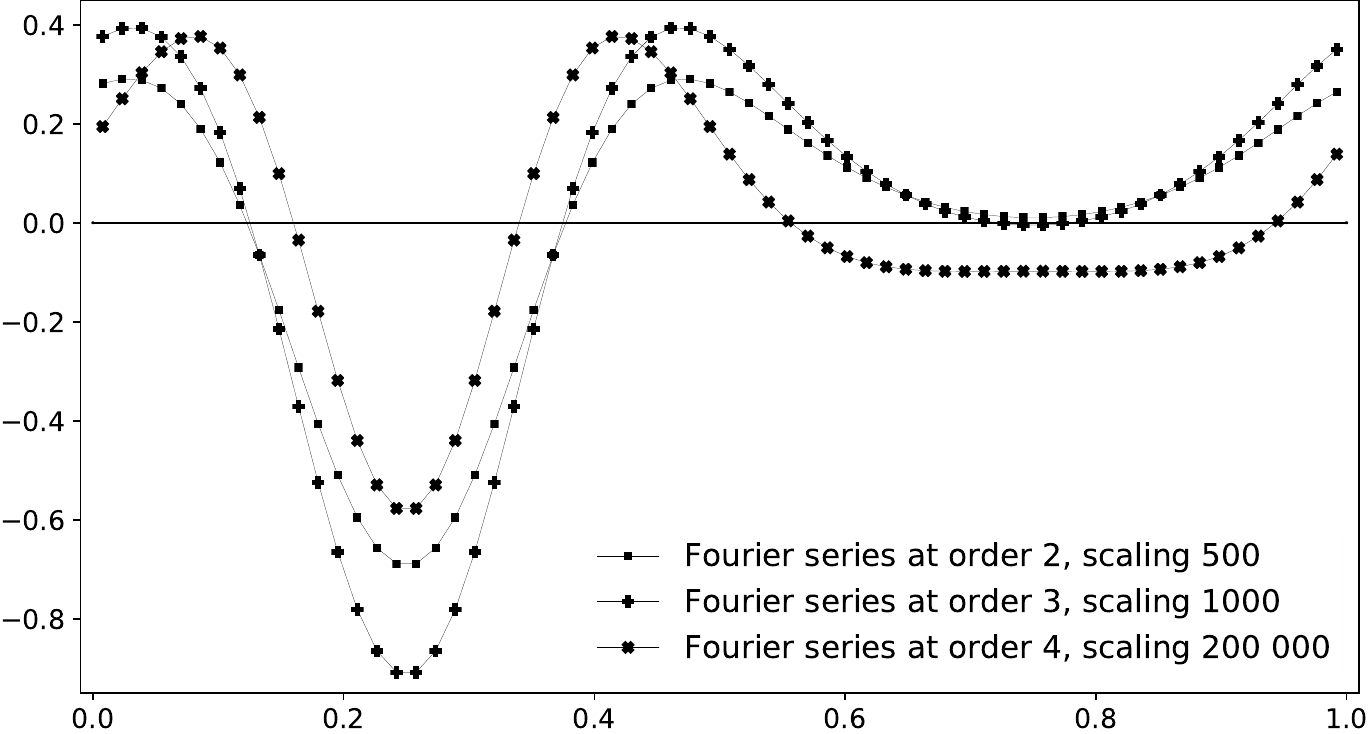}}
\caption{Stationary field, $ \, U = 0.005  $, 64 mesh points.
  Errors obtained for various levels of approxiation.}  
\label{fig-04-b} \end{figure} 

\smallskip  \noindent 
When $ \, U = 0.05 $, we  use  30 Fourier modes for  64 and 128 mesh points
and  60  modes  with 256 and 512 mesh points.
The results are presented  
in Table~\ref{tab-03}.
As advection speed increases, convergence also becomes increasingly difficult.
Nevertheless, the speed of convergence towards the stationary state is not directly correlated with
the order of the underlying partial differential equation.
For example, the numerical order of convergence for the fourth-order
partial differential equation is only  1.46.

\smallskip \noindent 
We have no complete explanation for  
the low orders of convergence reported in Tables \ref{tab-01} to~\ref{tab-03}.
Observe that the complete discrete dynamics (\ref{iterations-schema}) is diverging, as reported in Figure~\ref{modes-Gamma}, 
whereas the conserved moment  defined as usual by (\ref{densite-pierre}) is converging towards some given function
as the number of iterations tends to infinity.
Nevertheless, in studies of   
the finite difference method, it is well known (see {\it e.g.} \cite{LR56})
that the truncation error $ \, \theta_{\Delta t} \, $ is not identical to the error
$ \, \varepsilon_{\Delta t} $. We have in general for a global time $ \, T = N \Delta t \, $
composed of  $ \, N \, $ iterations,
\moneqstar 
\varepsilon_{\Delta t}  \leq C(T) \,\,  \theta_{\Delta t}  .
\monendstar 
  When the global integration time $ \, T \, $ is fixed (see Sections 6 to 8), we have a classical stability problem  
  and the error and the truncation errors have the same order of convergence.
This kind of analysis can  probably be extended for the high-order asymptotic expansion studied in this contribution.

\smallskip \noindent
When simulation time tends towards infinity, 
the discrete  dynamical 
system  (\ref{iterations-schema}) can  diverge exponentially as shown in Figure  \ref{modes-Gamma}.
In particular, some eigenmodes of the pure stationary lattice Boltzmann scheme correspond to an  unstable 
dynamics.
Nevertheless, the conserved quantity
defined in  (\ref{densite-pierre}) is converging towards a constant state.
The choice of this moment plays the role of a filter inside a diverging process.
As a result, the error analysis is difficult in the stationary  case. 
We can  reasonably assume that the  error $ \, \varepsilon_{\Delta t} \, $ and 
the   truncation error $ \, \theta_{\Delta t} \, $ behave as 
$ \,  \varepsilon_{\Delta t}  \simeq  C(T) \,\,  \theta_{\Delta t}  $.
Moreover,  the coefficient $ \, C(T) \, $ is in general tending towards infinity as $ \, T \, $ tends to 
infinity.
This kind of discrete stability violating the continuous stability criterion  
has been reported in the analysis of spectral methods
(see {\it e.g.} \cite{BM89}). Then the convergence accuracy is the result of the confrontation of a diverging
stability coefficient  and a converging truncation error. 
These remarks are not completely satisfactory. Nevertheless they offer the  
beginnings of an 
explanation for the curious convergence of stationary fields
compared to the convergence results for a finite-time evolution. 
  
\bigskip  \bigskip
\noindent {\bf \large    10) \quad  Conclusion} 

\smallskip \noindent 
In this contribution, we have extended the ABCD asymptotic analysis
developed in \cite{Du22}  and~\cite{DBL23} to a   inhomogeneous linear problem.
It has been necessary to develop a  library of Fourier series to
approximate with high accuracy the equivalent partial 
differential equations at orders 1 to 4.
The differential  operators have been explicated with the help of formal calculation
and, in particular, the Sagemath \cite{sage} library.

\smallskip \noindent
We have proven by numerical experiment  
that the asymptotic equivalent partial
differential equations constitute  a good approximation of the D1Q3 lattice Boltzmann scheme.
A major surprise in our study  concerns the  
finite-time evolution. We have put in evidence 
the importance of a correct initialization order 
to force the Boltzmann scheme to simulate 
a partial differential equation at high-order and obtain a convergence order
consistent with the formal approximation order.

\smallskip \noindent
For a  stationary problem after a long time evolution,
we put in evidence that  pure stationary modes of the lattice Boltzmann scheme 
can be unstable with a sinusoidal advective velocity. 
Nevertheless, 
the asymptotic expansion conducts to 
more and more precise approximations as the size of the mesh tends to zero.
The models suggested by the partial differential equations
are asymptotically correct but the order of accuracy is not the one
suggested by the order of the partial differential equation.
We leave this question to be addressed  
in a future work.

\smallskip \noindent
It would be also useful to consider a two-dimensional situation to get information
about anisotropic defects of lattice Boltzmann schemes.

%
\bigskip   \bigskip
\noindent {\bf  \large  Acknowledgments }

\smallskip 
This work has been supported by a public grant from the Fondation
Math\'ematique Jacques Hadamard as part of the ``Investissement
d'avenir'' project, reference ANR-11-LABX-0056-LMH, LabEx LMH.
FD thanks Thomas Bellotti for  an enlightening conversation on the order of initialization  of a multistep scheme.
FD  thanks also the Beijing Computational Science Research Center
and in particular  Li-Shi Luo 
for their  hospitality during the summer 2023. 
Last but not least, the referees proposed excellent suggestions that were not present
 in our initial text. These have been included in the final version of this work.

\bigskip   \bigskip
\noindent {\bf \large     Annex A. Proof of Proposition 3} 

\smallskip \noindent
We first recall the general result presented in \cite{Du22} and \cite{DBL23}.
We have 
\moneq \label{formules} \left \{ \begin {array}{rl}
\Gamma_1 (W)  & \!\!\! =     A \, W + B \, \Phi(W) \\
\Psi_1 (W) & \!\!\! =   \dd \Phi(W) .  \Gamma_1 (W)  - \big( C \, W + D \, \Phi(W)  \big) \\ 
\Gamma_2 (W) & \!\!\! =  B \, \Sigma \, \Psi_1 (W) \\
\Psi_2 (W) & \!\!\! =   \Sigma \,  \dd \Psi_1 (W) .   \Gamma_1 (W) + \dd \Phi(W) .  \Gamma_2 (W) 
- D \, \Sigma \,  \Psi_1 (W)   \\ 
\Gamma_3(W)   & \!\!\! = B \, \Sigma  \, \Psi_2 (W) + {{1}\over{12}}  B_2 \, \Psi_1  (W)
-  {{1}\over{6}} \, B \,  \dd \Psi_1 (W) .  \Gamma_1 (W) \\
\Psi_3 (W) & \!\!\! = \Sigma \, \dd \Psi_1 (W) .  \Gamma_2 (W)  +  \dd \Phi(W) .  \Gamma_3(W) -  D \, \Sigma \, \Psi_2 (W)
+ \Sigma \, \dd \Psi_2 (W) .  \Gamma_1 (W) \\
&  +{1\over6} \, D \, \dd \Psi_1 (W) .  \Gamma_1 (W)  - {1\over12} \, D_2 \, \Psi_1 (W)
- {1\over12} \, \dd \big( \dd \Psi_1 (W) .  \Gamma_1 (W) \big) . \Gamma_1 (W)  \\ 
\Gamma_4(W)  & \!\!\! =  B \, \Sigma \, \Psi_3 (W) + {1\over4} \, B_2 \, \Psi_2 (W)  +  {1\over6} \, B \, D_2 \, \Sigma \, \Psi_1 (W) 
 -   {1\over6} \, A \, B \, \Psi_2 (W) \\ & 
 -  {1\over6} \, B \, \dd \, (\dd \Phi . \Gamma_1 )  . \Gamma_2 (W)
 -  {1\over6} \, B \, \dd \, (\dd \Phi . \Gamma_2 ) . \Gamma_1 (W)  \\ &  
 - {1\over6} \, B \, \Sigma \, \dd \, (\dd \Psi_1 (W) .  \Gamma_1 ) .  \Gamma_1 (W) . 
\end {array} \right. \monend
With the one-dimensional relations (\ref{Abarre-et-al}), we have in particular
\moneqstar
B \, \Phi =  {\overline B} \, \partial_x   \Phi  = {\overline B} \,\,  \partial_x (E(x) \, W) =  {\overline B} \, \delta \, W .
\monendstar
Then $ \,\, \Gamma_1 = {\overline A} \, \partial_x  W +  {\overline B} \, \delta \, W = \alpha_1 \, W \,\, $
with $ \,\, \alpha_1 = {\overline A} \,\, \partial_x + {\overline B} \,\, \delta \,\, $ and the first relation of
the family (\ref{operateurs-barres}) is proven.
We have as previously $ \, D \, \Phi =  {\overline D} \,\,  \partial_x (E(x) \, W) =   {\overline D} \, \delta \, W \, $ and 

\smallskip \noindent $
\Psi_1 (W) =   \dd \Phi(W) .  \Gamma_1 (W)  - \big( C \, W + D \, \Phi(W)  \big) = E \, \alpha_1 \,   W
- ( {\overline C} \, \partial_x W +  {\overline D} \,  \delta \, W ) $

\smallskip \noindent $ \quad \qquad \,\, 
= \big[ E \, \alpha_1 - ( {\overline C} \,\, \partial_x + {\overline D} \,\, \delta ) \big] \, W \equiv \beta_1 \, W $.

\smallskip \noindent
Then  the second relation in (\ref{operateurs-barres}) relative to $ \, \beta_1 \, $ is established.

\smallskip \noindent
From the fact that the matrix $ \, \Sigma \, $ is constant, we have now 

\smallskip \noindent
$ \, \Gamma_2 (W) =  B \, \Sigma \, \Psi_1 (W) =  {\overline B} \, \partial_x \, \Sigma \, \beta_1 \, W = 
           {\overline B} \, \Sigma \, \partial_x \,  \beta_1 \, W \equiv \alpha_2 \, W \, $
           
\smallskip \noindent
and the relation in   (\ref{operateurs-barres}) relative to $ \, \alpha_2 \, $ is proven.

\smallskip \noindent 
When we differentiate the vector field $ \, \Psi_1 (W) $, we have
$ \,\, \dd  \Psi_1 (W) . \xi = \beta_1 \, \xi \,\, $ and

\noindent 
$ \, \dd \Psi_1 (W) .   \Gamma_1 (W) = \beta_1 \, \alpha_1 \, W $. Then

\smallskip \noindent $
\Psi_2 (W) =   \Sigma \,  \dd \Psi_1 (W) .   \Gamma_1 (W) + \dd \Phi(W) .  \Gamma_2 (W) - D \, \Sigma \,  \Psi_1 (W) $ 

\smallskip \noindent $\quad \qquad \,\, =
\Sigma \, \beta_1 \, \alpha_1 \, W + E \, \alpha_2 \, W - {\overline D} \, \delta_x \, \Sigma \, \beta_1 \, W
= \big( \Sigma \, \beta_1 \, \alpha_1  + E \, \alpha_2  - {\overline D} \, \Sigma \, \partial_x \, \beta_1 \big) \, W \equiv \beta_2 \, W $

\smallskip \noindent
and the relation relation to $ \, \beta_2 \, $ in   (\ref{operateurs-barres}) is proven.
We have also

\smallskip \noindent $
B_2 = (A \, B + B\, D ) =  {\overline A} \,\, \partial_x \,\,  {\overline B} \,\, \partial_x  + 
{\overline B} \,\, \partial_x \,\,  {\overline D} \,\, \partial_x = \big(  {\overline A} \,\,  {\overline B} +  {\overline B} \,\,  {\overline D}
\big) \,  \partial_x^2 = {\overline B_2}  \,  \partial_x^2 $.

\smallskip \noindent
In consequence, 

\smallskip \noindent $
\Gamma_3(W)  = B \, \Sigma  \, \Psi_2 (W) + {{1}\over{12}}  B_2 \, \Psi_1  (W) -  {{1}\over{6}} \, B \,  \dd \Psi_1 (W) .  \Gamma_1 (W) $ 

\smallskip \noindent $\quad \qquad  =
 {\overline B} \, \Sigma \, \partial_x \, \beta_2 \, W  + {{1}\over{12}}  {\overline B_2}  \,  \partial_x^2 \, \beta_1 \, W
 -  {{1}\over{6}} \, {\overline B} \, \partial_x \, \beta_1 \, \alpha_1 \, W $

\smallskip \noindent $\quad \qquad  = 
\big(  {\overline B} \, \Sigma \, \partial_x \, \beta_2   + {{1}\over{12}}  {\overline B_2}  \,  \partial_x^2 \, \beta_1 
 -  {{1}\over{6}} \, {\overline B} \, \partial_x \, \beta_1 \, \alpha_1 \big) \, W  \equiv \alpha_3 \, W  $ 

 \smallskip \noindent
 and the expression of the operator $ \, \alpha_3 \, $ in   (\ref{operateurs-barres}) is established.
From the relation (\ref{formules}), we have 
$ \,  \dd \big( \dd \Psi_1 (W) .  \Gamma_1 (W) \big) . \xi =  \dd \big(  \beta_1 \, \alpha_1 \, W \big) . \, \xi \, $
and

\smallskip \noindent $
\dd \big( \dd \Psi_1 (W) .  \Gamma_1 (W) \big) . \Gamma_1 (W)  =  \beta_1 \, \alpha_1 \,  \alpha_1 \, W
= \beta_1 \, \alpha_1^2  \,  W $. Then

\smallskip \noindent $
\Psi_3 (W) = \Sigma \, \dd \Psi_1 (W) .  \Gamma_2 (W)  +  \dd \Phi(W) .  \Gamma_3(W) -  D \, \Sigma \, \Psi_2 (W)
+ \Sigma \, \dd \Psi_2 (W) .  \Gamma_1 (W) $

\smallskip \noindent $\qquad \qquad  +
{1\over6} \, D \, \dd \Psi_1 (W) .  \Gamma_1 (W)  - {1\over12} \, D_2 \, \Psi_1 (W)
- {1\over12} \, \dd \big( \dd \Psi_1 (W) .  \Gamma_1 (W) \big) . \Gamma_1 (W) $ 

\smallskip \noindent $\quad \qquad \, 
= \Sigma \, \beta_1 \, \alpha_2 \, W  + E \, \alpha_3 \, W -  {\overline D} \, \partial_x \, \Sigma \, \beta_2 \, W 
+  \Sigma \, \beta_2 \, \alpha_1 \, W $

\smallskip \noindent $\qquad \qquad 
+ {1\over6} \, {\overline D} \, \partial_x \, \beta_1 \, \alpha_1 \, W
  - {1\over12} \,  {\overline D_2} \, \partial_x^2 \, \beta_1 \, W  - {1\over12} \, \beta_1 \, \alpha_1^2  \,  W $

\smallskip \noindent $\quad \qquad  \, =
\big(   \Sigma \, \beta_1 \, \alpha_2 + E \, \alpha_3 -  {\overline D} \, \partial_x \, \Sigma \, \beta_2 
+  \Sigma \, \beta_2 \, \alpha_1 +  {1\over6} \, {\overline D} \, \partial_x \, \beta_1 \, \alpha_1 
  - {1\over12} \,  {\overline D_2} \, \partial_x^2 \, \beta_1  - {1\over12} \, \beta_1 \, \alpha_1^2 \big) \, W $ 
 
\smallskip \noindent
and the relation   (\ref{operateurs-barres}) concerning $ \, \beta_3 \, $ is established.

\smallskip \noindent
We observe now  that

\smallskip \noindent $
B \, \dd \, (\dd \Phi . \Gamma_2 ) . \, \xi =  {\overline B} \,\partial_x \, \dd ( E \, \alpha_2 \, W ) . \, \xi
=  {\overline B} \, \delta \, \alpha_2 \, \xi \,\, $ and
$ \, B \, \dd \, (\dd \Phi . \Gamma_2 ) . \Gamma_1 (W) =   {\overline B} \, \delta \, \alpha_2 \, \alpha_1 \, W $.

\smallskip \noindent 
We have finally

\smallskip \noindent $
\Gamma_4(W) =  B \, \Sigma \, \Psi_3 (W) + {1\over4} \, B_2 \, \Psi_2 (W)  +  {1\over6} \, B \, D_2 \, \Sigma \, \Psi_1 (W) 
-   {1\over6} \, A \, B \, \Psi_2 (W) $

\smallskip \noindent $\qquad \qquad
 -  {1\over6} \, B \, \dd \, (\dd \Phi . \Gamma_1 )  . \Gamma_2 (W)
 -  {1\over6} \, B \, \dd \, (\dd \Phi . \Gamma_2 ) . \Gamma_1 (W)
 - {1\over6} \, B \, \Sigma \, \dd \, (\dd \Psi_1 (W) .  \Gamma_1 ) .  \Gamma_1 (W) $

\smallskip \noindent $\quad \qquad  \, =
{\overline B} \,  \Sigma \, \partial_x \, \beta_3  \, W  +  {1\over4} \,  {\overline B_2} \, \partial_x^2 \, \beta_2 \, W 
+  {1\over6} \, {\overline B} \, {\overline D_2} \,  \Sigma \,\partial_x^3 \,  \beta_1 \, W 
-  {1\over6} \, {\overline A} \, {\overline B} \,\partial_x^2 \,  \beta_2 \, W $

\smallskip \noindent $\qquad \qquad
-  {1\over6} \,{\overline B} \, \delta \, \alpha_1 \, \alpha_2 \, W 
-  {1\over6} \,{\overline B} \, \delta \, \alpha_2 \, \alpha_1 \, W
-  {1\over6} \,{\overline B} \,  \Sigma \, \partial_x \, \beta_1 \, \alpha_1^2 \, W $ 

\smallskip \noindent $\quad \qquad  \, =
\Big[ {\overline B} \,  \Sigma \, \partial_x \, \beta_3   +  {1\over4} \,  {\overline B_2} \, \partial_x^2 \, \beta_2 
+  {1\over6} \, {\overline B} \, {\overline D_2} \,  \Sigma \,\partial_x^3 \,  \beta_1 
-  {1\over6} \, {\overline A} \, {\overline B} \,\partial_x^2 \,  \beta_2  
-  {1\over6} \,{\overline B} \, \delta \, \alpha_1 \, \alpha_2 $ 

\smallskip \noindent $\qquad \qquad
-  {1\over6} \,{\overline B} \, \delta \, \alpha_2 \, \alpha_1
-  {1\over6} \,{\overline B} \,  \Sigma \, \partial_x \, \beta_1 \, \alpha_1^2 \Big] \, W
\equiv \alpha_4 \, W $

\smallskip \noindent
and the last relation of (\ref{operateurs-barres}) 
giving the operator $ \, \alpha_4 \, $ is explicated.
\hfill $\square $

%

\newpage 
\bigskip \bigskip    \noindent {\bf \large     Annex B. Proof of Proposition 4} 

\smallskip \noindent
We start from the relations  (\ref{operateurs-barres}) that we write again which we rewrite for clarity of reading:
\moneqstar  \left\{ \begin{array} {l}  
  \alpha_1 =  {\overline A} \,\, \partial_x + {\overline B} \,\, \delta \\
  \beta_1 = E \, \alpha_1 - ( {\overline C} \,\, \partial_x + {\overline D} \,\, \delta ) \\
  \alpha_2 =  {\overline B} \, \Sigma \, \partial_x \, \beta_1  \\
  \beta_2 = \Sigma \, \beta_1 \, \alpha_1 + E \, \alpha_2 - {\overline D} \,  \Sigma \, \partial_x \, \beta_1 \\ 
   \alpha_3 =  {\overline B} \,  \Sigma \, \partial_x \, \beta_2 
   + {1\over12} \,  {\overline {B_2}} \, \partial_x^2 \,  \beta_1  -  {1\over6} \,  {\overline B} \, \partial_x \, \beta_1 \, \alpha_1 \\
 \beta_3 = \Sigma \, \beta_1 \, \alpha_2 + E \, \alpha_3 -  {\overline D} \,  \Sigma \, \partial_x \, \beta_2
 +  \Sigma \, \beta_2 \, \alpha_1 + {1\over6} \, {\overline D} \,   \partial_x \, \beta_1 \, \alpha_1
 -{1\over12} \,  \beta_1 \, \alpha_1^2  -{1\over12} \,  {\overline D_2} \,   \partial_x^2 \, \beta_1 \\  
 \alpha_4 =  {\overline B} \,  \Sigma \, \partial_x \, \beta_3 +  {1\over4} \,  {\overline B_2} \, \partial_x^2 \, \beta_2
+  {1\over6} \, {\overline B} \, {\overline D_2} \,  \Sigma \,\partial_x^3 \,  \beta_1
-  {1\over6} \, {\overline A} \, {\overline B} \,\partial_x^2 \,  \beta_2 \\ \qquad 
-  {1\over6} \,{\overline B} \, \delta \, \alpha_1 \, \alpha_2
-  {1\over6} \,{\overline B} \, \delta \, \alpha_2 \, \alpha_1  -  {1\over6} \,{\overline B} \,  \Sigma \, \partial_x \, \beta_1 \, \alpha_1^2  \,. 
\end{array} \right. \monendstar 
With the D1Q3 lattice Boltzmann scheme, we have 
\moneqstar 
\delta =  \begin{pmatrix} \lambda  \, \partial_u \\ \lambda^2 \, \alpha \, \partial_x \end{pmatrix} \,,\,\,  
E(x) = \begin{pmatrix} \lambda \, U \, \cos(k\,x) \\ \lambda^2 \, \alpha \end{pmatrix} \,,\,\,
\Sigma =  \begin{pmatrix} \sigma & 0 \\ 0 &  \sigma'  \end{pmatrix} 
\monendstar
and
\moneqstar
{\overline A} = 0 \,,\,\,  {\overline B} = \big( 1 \,,\, 0 \big)  \,,\,\,
{\overline C} = \begin{pmatrix} {{2\, \lambda^2}\over{3}} \\ 0 \end{pmatrix}  \,,\,\,
{\overline D} = \begin{pmatrix} 0 & {1\over3} \\ \lambda^2 & 0  \end{pmatrix} .
\monendstar 
Then we obtain 
$ \, \alpha_1 = \delta = \lambda \, \partial_u $,

\smallskip \noindent
$\, \beta_1 =  \begin{pmatrix}  u \\ \lambda^2 \, \alpha \end{pmatrix}  \lambda \, \partial_u  
- \begin{pmatrix} {{2\, \lambda^2}\over{3}} \\ 0 \end{pmatrix} \partial_x
- \begin{pmatrix} 0 & {1\over3} \\ \lambda^2 & 0  \end{pmatrix}  \begin{pmatrix} \lambda  \, \partial_u \\ \lambda^2 \, \alpha \, \partial_x \end{pmatrix}
=  \begin{pmatrix}  \lambda \, u \,   \partial_u - {{2}\over{3}} \,  \lambda^2 \, \partial_x -  {{\lambda^2}\over{3}} \,\alpha \, \partial_x \\
 \lambda^3 \, \alpha \,  \partial_u  -  \lambda^3 \,   \partial_u \end{pmatrix} $

\smallskip \noindent $\quad \,\, 
=  \begin{pmatrix} \lambda  \, u \,  \partial_u -  {{1}\over{3}} \, \lambda^2 \, (\alpha+2) \, \partial_x \big) \\
 \lambda^3  \, (\alpha - 1) \,  \partial_u  \end{pmatrix} $

\smallskip \noindent
and the relation (\ref{alpha1-beta1}) is proven. We have for  second order accuracy
$\,  {\overline B} \,\, \Sigma = \big( \sigma  ,\, 0 \big) $. Then 

\smallskip \noindent
$ \alpha_2 = \sigma \,  \partial_x \, \big[ \lambda \, u \,  \partial_u -  {{\lambda^2}\over{3}} \, (\alpha+2) \, \partial_x \big) \big] =
\lambda^2  \, \sigma \,  \big( \partial_u^2 -  {{\alpha+2}\over{3}} \,  \partial_x^2 \big) $. 

\smallskip \noindent
For the microscopic variables, we have

\smallskip \noindent $
\beta_2 = \begin{pmatrix} \sigma & 0 \\ 0 &  \sigma'  \end{pmatrix} \lambda^2
\begin{pmatrix}  \lambda \, u \,   \partial_u - {{2}\over{3}} \,  \lambda^2 \, \partial_x -  {{\lambda^2}\over{3}} \,\alpha \, \partial_x \\
  \lambda^3 \, \alpha \,  \partial_u  -  \lambda^3 \,   \partial_u \end{pmatrix} \,  \lambda \, \partial_u
+ \begin{pmatrix}  u \\ \lambda^2 \, \alpha \end{pmatrix}  \lambda^2 \, \sigma \, \big( \partial_u^2 -  {{\alpha+2}\over{3}} \,  \partial_x^2 \big) $

\smallskip \noindent \qquad $
-  \lambda \, \begin{pmatrix} 0 & {{\sigma'}\over{3}} \\ \lambda^2 \, \sigma & 0 \end{pmatrix} \, \begin{pmatrix} 
u \,  \partial_u -  {{1}\over{3}} \, \lambda^2 \, (\alpha+2) \, \partial_x \\  \lambda^2  \, (\alpha - 1) \,  \partial_u \end{pmatrix} . $

\smallskip \noindent 
For the first component,

\smallskip \noindent $
\beta_2^1 = \sigma \, \big(  \lambda \, u \,   \partial_u - {{\alpha+2}\over{3}} \,  \lambda^2 \, \partial_x \big) \, \lambda \, \partial_u 
+  \lambda^2 \, u \, \sigma \, \big(\partial_u^2 -  {{\alpha+2}\over{3}} \,  \partial_x^2 \big)
-  \lambda^3 \,  {{\sigma'}\over{3}} \,  \partial_x  \,  (\alpha - 1) \,  \partial_u   $

\smallskip \noindent $ \quad \,\, 
= 2 \, \lambda^2 \, \sigma \,  u \,  \partial_u^2 -  {{\alpha+2}\over{3}} \,  \lambda^3 \, \partial_x \,  \partial_u
- \lambda^2 \,  {{\alpha+2}\over{3}} \, u \, \partial_x^2 - \lambda^3 \, (\alpha-1) \,  {{\sigma'}\over{3}} \,  \partial_x \,  \partial_u $ 

\smallskip \noindent $ \quad \,\, 
=  \lambda^3 \, \big[  2 \, \sigma \, {{u}\over{\lambda}} \,  \partial_u^2
- \big(  {{\alpha+2}\over{3}} \, \sigma -  {{\alpha-1}\over{3}} \,\sigma' \big) \,  \partial_x \,  \partial_u
-    {{\alpha+2}\over{3}} \,  \sigma \,  {{u}\over{\lambda}} \, \partial_x^2 \big] $ 

\smallskip \noindent 
and for the second component 

\smallskip \noindent $
\beta_2^2 =  \lambda^2 \, \sigma' \, \big( \lambda \, (\alpha-1) \,  \partial_u \big) \, \lambda \,  \partial_u
+ \lambda^4 \, \alpha \, \sigma \,  \big( \partial_u^2 -  {{\alpha+2}\over{3}} \,  \partial_x^2 \big)
- \lambda^4 \, \sigma \,  \big(  {{u}\over{\lambda}} \, \partial_u -  {{\alpha+2}\over{3}} \, \partial_x \big) $

\smallskip \noindent $ \quad \,\, 
= \lambda^4 \, \big[  (\alpha-1) \, \sigma' \, \partial_u^2  +  \alpha \, \sigma \,  \big( \partial_u^2 -  {{\alpha+2}\over{3}} \,  \partial_x^2 \big) 
- \sigma \,   \big( \partial_u^2 -  {{\alpha+2}\over{3}} \, \partial_x^2 \big) \big] $

\smallskip \noindent $ \quad \,\, 
= \lambda^4 \,   (\alpha-1) \, \big(  (\sigma + \sigma') \, \partial_u^2 -  {{\alpha+2}\over{3}} \,\sigma \, \partial_x^2 \big) $.

\smallskip \noindent 
Then the relations  (\ref{alpha2-beta2}) are established. At third order, we have from  (\ref{operateurs-barres}),

\smallskip \noindent $
\alpha_3 =  {\overline B} \,  \Sigma \, \partial_x \, \beta_2 
   + {1\over12} \,  {\overline {B_2}} \, \partial_x^2 \,  \beta_1  -  {1\over6} \,  {\overline B} \, \partial_x \, \beta_1 \, \alpha_1  $ 

\smallskip \noindent $ \quad \,\, 
= \sigma \, \partial_x \beta_2^1    + {1\over12} \, \big( 0  ,\,  {{1}\over{3}} \big) \, \partial_x \beta_1
-  {1\over6} \,  {\overline B} \, \partial_x \begin{pmatrix}  \lambda \, u \,   \partial_u - {{2}\over{3}} \,  \lambda^2 \, \partial_x
-  {{\lambda^2}\over{3}} \,\alpha \, \partial_x \\  \lambda^3 \, \alpha \,  \partial_u
-  \lambda^3 \,   \partial_u \end{pmatrix} \, \lambda \, \partial_u $ 

\smallskip \noindent $ \quad \,\, 
=    \lambda^3 \,    \sigma \, \partial_x \big[  2 \, \sigma \, {{u}\over{\lambda}} \,  \partial_u^2
- \big(  {{\alpha+2}\over{3}} \, \sigma -  {{\alpha-1}\over{3}} \,\sigma' \big) \,  \partial_x \,  \partial_u
-    {{\alpha+2}\over{3}} \,  \sigma \,  {{u}\over{\lambda}} \, \partial_x^2 \big]
  + {1\over12} \,  {{\lambda^3}\over{3}} \, (\alpha-1) \, \partial_x^2 \, \partial_u  $

\smallskip \noindent $ \quad \quad 
-    {{\lambda}\over{6}} \, \partial_x \big[  \lambda \, u \,   \partial_u - {{\alpha+2}\over{3}} \,  \lambda^2 \, \partial_x \big] 
\, \lambda \,  \partial_u  $

\smallskip \noindent $ \quad \,\, 
=  \lambda^3  \Big[ \Big( 2 \, \sigma^2 - {1\over6} \Big) \, \partial_u^3  +
  \big[ -\sigma \, \big( {{\alpha+2}\over{3}} \, \sigma +  {{\alpha-1}\over{3}} \, \sigma' \big)  +  {{\alpha-1}\over{36}}
    + {{\alpha+2}\over{18}} \big] \,  \partial_x^2 \, \partial_u  
-  {{\alpha+2}\over{3}} \, \sigma^2 \, \partial_u  \,  \partial_x^2 $ 

\smallskip \noindent $ \quad \,\, 
= \lambda^3  \Big[ \Big( 2 \, \sigma^2 - {1\over6} \Big) \, \partial_u^3 
  + \Big(  {{\alpha+2}\over{3}}\, \big( {1\over6}-\sigma^2 \big) + {{\alpha-1}\over{3}}\,  \big( {1\over12}- \sigma\, \sigma' \big)
  \Big) \,  \partial_x^2 \,\partial_u    - {{\alpha+2}\over{3}}\,\sigma^2 \, \partial_u \, \partial_x^2   \Big] $

\smallskip \noindent
and the relation (\ref{alpha3}) is proven. We have now

\smallskip \noindent $ 
\beta_3 = \Sigma \, \beta_1 \, \alpha_2 + E \, \alpha_3 -  {\overline D} \,  \Sigma \, \partial_x \, \beta_2
 +  \Sigma \, \beta_2 \, \alpha_1 + {1\over6} \, {\overline D} \,   \partial_x \, \beta_1 \, \alpha_1
 -{1\over12} \,  \beta_1 \, \alpha_1^2  - {1\over12} \,  {\overline D_2} \,   \partial_x^2 \, \beta_1  $ 

 \smallskip \noindent
and we can precise these seven terms: 

\smallskip \noindent $ 
\Sigma \, \beta_1 \, \alpha_2 =   \begin{pmatrix} \sigma & 0 \\ 0 &  \sigma'  \end{pmatrix} \lambda \,
 \begin{pmatrix}  \lambda^3 \, \big[  2 \, \sigma \, {{u}\over{\lambda}} \,  \partial_u^2
- \big(  {{\alpha+2}\over{3}} \, \sigma -  {{\alpha-1}\over{3}} \,\sigma' \big) \,  \partial_x \,  \partial_u
-    {{\alpha+2}\over{3}} \,  \sigma \,  {{u}\over{\lambda}} \, \partial_x^2 \big] \\
   \lambda^4 \,   (\alpha-1) \, \big(  (\sigma + \sigma') \, \partial_u^2 -  {{\alpha+2}\over{3}} \,\sigma \, \partial_x^2 \big) \end{pmatrix}
 \lambda^2 \, \sigma \,  \Big(  \partial_u^2 -  {{\alpha+2}\over{3}} \, \partial_x^2 \Big)   $

\smallskip \noindent $ 
E \, \alpha_3 = \begin{pmatrix} u \\ \lambda^2 \, \alpha \end{pmatrix} \, 
 \lambda^3  \Big[ \Big( 2 \, \sigma^2 - {1\over6} \Big) \, \partial_u^3 
  + \Big(  {{\alpha+2}\over{3}}\, \big( {1\over6}-\sigma^2 \big) + {{\alpha-1}\over{3}}\,  \big( {1\over12}- \sigma\, \sigma' \big)
  \Big) \,  \partial_x^2 \,\partial_u    - {{\alpha+2}\over{3}}\,\sigma^2 \, \partial_u \, \partial_x^2   \Big]  $

\smallskip \noindent $ 
-  {\overline D} \,  \Sigma \, \partial_x \, \beta_2 = - \begin{pmatrix} 0 & {1\over3} \\ \lambda^2 & 0  \end{pmatrix}
\begin{pmatrix} \sigma & 0 \\ 0 & \sigma' \end{pmatrix} \partial_x
\begin{pmatrix}  \lambda^3 \, \big[  2 \, \sigma \, {{u}\over{\lambda}} \,  \partial_u^2
- \big(  {{\alpha+2}\over{3}} \, \sigma -  {{\alpha-1}\over{3}} \,\sigma' \big) \,  \partial_x \,  \partial_u
-    {{\alpha+2}\over{3}} \,  \sigma \,  {{u}\over{\lambda}} \, \partial_x^2 \big] \\ 
\lambda^4 \,   (\alpha-1) \, \big(  (\sigma + \sigma') \, \partial_u^2 -  {{\alpha+2}\over{3}} \,\sigma \, \partial_x^2 \big)  \end{pmatrix} $

\smallskip \noindent $ 
\Sigma \, \beta_2 \, \alpha_1 = \begin{pmatrix} \sigma & 0 \\ 0 & \sigma' \end{pmatrix} \,
\begin{pmatrix}  \lambda^3 \, \big[  2 \, \sigma \, {{u}\over{\lambda}} \,  \partial_u^2
- \big(  {{\alpha+2}\over{3}} \, \sigma -  {{\alpha-1}\over{3}} \,\sigma' \big) \,  \partial_x \,  \partial_u
-    {{\alpha+2}\over{3}} \,  \sigma \,  {{u}\over{\lambda}} \, \partial_x^2 \big] \\ 
  \lambda^4 \,   (\alpha-1) \, \big(  (\sigma + \sigma') \, \partial_u^2 -  {{\alpha+2}\over{3}} \,\sigma \, \partial_x^2 \big)  \end{pmatrix}
\, \lambda \, \partial_u $ 

\smallskip \noindent $ 
{1\over6} \, {\overline D} \,   \partial_x \, \beta_1 \, \alpha_1 = {1\over6} \, \begin{pmatrix} 0 & {1\over3} \\ \lambda^2 & 0  \end{pmatrix}
\,   \partial_x \, \begin{pmatrix} \lambda  \, u \,  \partial_u -  {{1}\over{3}} \, \lambda^2 \, (\alpha+2) \, \partial_x \big) \\
 \lambda^3  \, (\alpha - 1) \,  \partial_u  \end{pmatrix}\, \lambda \, \partial_u $ 

\smallskip \noindent $ 
-{1\over12} \,  \beta_1 \, \alpha_1^2 = -{1\over12} \,
\begin{pmatrix} \lambda  \, u \,  \partial_u -  {{1}\over{3}} \, \lambda^2 \, (\alpha+2) \, \partial_x \big) \\
 \lambda^3  \, (\alpha - 1) \,  \partial_u  \end{pmatrix} \, \lambda^2 \, \partial_u^2 $ 

\smallskip \noindent $ 
- {1\over12} \,  {\overline D_2} \,   \partial_x^2 \, \beta_1 =  - {{\lambda^2}\over12} \, \begin{pmatrix} 1 & 0 \\ 0 & {1\over3} \end{pmatrix} \, 
\partial_x^2 \, \begin{pmatrix} \lambda  \, u \,  \partial_u -  {{1}\over{3}} \, \lambda^2 \, (\alpha+2) \, \partial_x \big) \\
 \lambda^3  \, (\alpha - 1) \,  \partial_u  \end{pmatrix}  . $

\smallskip \noindent
Then the first component of $ \, \beta_3 \, $ is given by the relation

\smallskip \noindent $ 
\beta_3^1 = \lambda^4 \, \Big[ {{\alpha+2}\over9}\, \big[ -(1-\alpha)\, \sigma \, \sigma'  +  (\big (\alpha +2)\, \sigma^2  + {1\over4} \big) \big] \, \partial_x^3   $

\smallskip \noindent $  \quad \quad
+ \, U \,  \big[  -2 \,  {{\alpha+2}\over3}\, \sigma^2 + {{1-\alpha}\over3} \,\sigma \, \sigma' +   {{1+\alpha}\over12} \big] 
\, \partial_x^2  \, \partial_u - 2\, U \, {{\alpha+2}\over3} \, \sigma^2 \, \partial_u \, \partial_x^2 $

\smallskip \noindent $  \quad \quad 
 +  \big[ -2\,  {{\alpha+2}\over3}\, \sigma^2
  +  {{1-\alpha}\over3} \,(2 \, \sigma \, \sigma' + \sigma'^2 -{1\over4}) \big] \, \partial_x \, \partial_u^2 
 + \big( 5 \, \sigma^2  -{1\over4} \big) \, U \,  \partial_u^3 \Big] $

\smallskip \noindent
and the second is given by 

\smallskip \noindent $ 
\beta_3^2 = \lambda^5 \, \Big[ {{1-\alpha}\over3} \,  \big[  (\alpha+2) \, \sigma^2 + (1 + 2\, \alpha) \, \sigma \, \sigma'
  -  {{1+\alpha}\over4} \big] \, \partial_x^2 \, \partial_u 
+ (1-\alpha) \,  {{\alpha+2}\over3}\, \sigma \, ( \sigma +  \sigma' ) \, \partial_u \, \partial_x^2 $ 

\smallskip \noindent $  \quad \quad 
- (1-\alpha) \, \big( 2 \,  \sigma^2 + 2 \, \sigma \, \sigma' + \sigma'^2  -{1\over4}  \big) \, \partial_u^3  \Big]  \, . $

\smallskip \noindent
The relation (\ref{beta3}) is proven. Finally,

\smallskip \noindent $
\alpha_4 =  {\overline B} \,  \Sigma \, \partial_x \, \beta_3 +  {1\over4} \,  {\overline B_2} \, \partial_x^2 \, \beta_2
+  {1\over6} \, {\overline B} \, {\overline D_2} \,  \Sigma \,\partial_x^3 \,  \beta_1
-  {1\over6} \, {\overline A} \, {\overline B} \,\partial_x^2 \,  \beta_2 
-  {1\over6} \,{\overline B} \, \delta \, \alpha_1 \, \alpha_2 $

\smallskip \noindent $  \quad \quad 
-  {1\over6} \,{\overline B} \, \delta \, \alpha_2 \, \alpha_1  -  {1\over6} \,{\overline B} \,  \Sigma \, \partial_x \, \beta_1 \, \alpha_1^2  \,. $

\smallskip \noindent
After some  lines of algebra, 

\smallskip \noindent $ 
{{1}\over{\lambda^4}} \, \alpha_4 =  \big[ \big( {{\alpha+2}\over3} \big)^2 \, \sigma^3 + (\alpha-1)\, {{\alpha+2}\over9} \, \sigma^2 \, \sigma'
- {{\alpha+2}\over36} \, \alpha \, \sigma \big] \, \partial_x^4 $

\smallskip \noindent $  \qquad \quad 
+  \big[  -2\, {{\alpha+2}\over3} \, \sigma^3 + 2 \, {{1-\alpha}\over3} \, \sigma^2 \, \sigma' +  {{1-\alpha}\over3} \, \sigma \, \sigma'^2
+  {{1+ 2\,\alpha}\over9} \, \sigma + {{\alpha-1}\over12} \, \sigma' \big] \, \partial_x^2 \, \partial_u^2 $  

\smallskip \noindent $  \qquad \quad 
+  \big[  -2\, {{\alpha+2}\over3} \, \sigma^3 +  {{1-\alpha}\over3} \, \sigma^2 \, \sigma' +  {{7 + 5\,\alpha}\over36} \, \sigma
 \big] \, \partial_u \, \partial_x^2 \, \partial_u $ $  
+  \, {{\alpha+2}\over3} \, \sigma \, \big( -2 \, \sigma + {1\over6} \big) \, \partial_u^2 \, \partial_x^2
+ \sigma \, \big( 5 \, \sigma^2 - {3\over4} \big) \,  \partial_u^4   $

\smallskip \noindent $  \quad \quad \,\, 
=  \big[ {{\alpha+2}\over9} \, \big( (\alpha+2) \, \sigma^3 - (1-\alpha)\, \sigma^2 \, \sigma'
  - {{\alpha}\over4} \, \sigma \big) \,  \partial_x^4   $
  
\smallskip \noindent $  \qquad \quad 
+ \,  \big[ -2\, {{\alpha+2}\over3} \, \sigma^3
+  {{1-\alpha}\over3} \, (2 \, \sigma^2 \, \sigma' + \sigma \, \sigma'^2 - {1\over4} \, \sigma' )
+ {{1 + 2\, \alpha}\over9} \,   \sigma \big] \,  \partial_x^2 \, \partial_u^2  $

\smallskip \noindent $  \qquad \quad 
+  \big[  -2\, {{\alpha+2}\over3} \, \sigma^3 +  {{1-\alpha}\over3} \, \sigma^2 \, \sigma' +  {{7 + 5\,\alpha}\over36} \, \sigma
 \big] \, \partial_u \, \partial_x^2 \, \partial_u $ $  
+  \, {{\alpha+2}\over3} \, \sigma \, \big( -2 \, \sigma + {1\over6} \big) \, \partial_u^2 \, \partial_x^2
+ \sigma \, \big( 5 \, \sigma^2 - {3\over4} \big) \,  \partial_u^4   $

\smallskip \noindent
and the relation (\ref{alpha4}) is established.
This completes the proof. \hfill $\square$

\bigskip \bigskip      \noindent {\bf  \large  References }


\end{document}